\documentclass[12pt]{amsart}
\usepackage[utf8]{inputenc}
\usepackage[all]{xy}
\usepackage{hyperref, amsmath, amssymb, amsthm, mathtools, fullpage, enumitem,mathrsfs, aliascnt}

\theoremstyle{plain}
\newtheorem{theorem}{Theorem}[section]
\newaliascnt{prop}{theorem}
\newtheorem{prop}[prop]{Proposition}
\aliascntresetthe{prop}

\newaliascnt{cor}{theorem}
\newtheorem{cor}[cor]{Corollary}
\aliascntresetthe{cor}

\newaliascnt{lemma}{theorem}
\newtheorem{lemma}[lemma]{Lemma}
\aliascntresetthe{lemma}

\newtheorem{mainthm}{Theorem}

\theoremstyle{definition}
\newaliascnt{defn}{theorem}
\newtheorem{defn}[defn]{Definition}
\aliascntresetthe{defn}

\newaliascnt{rem}{theorem}
\newtheorem{rem}[rem]{Remark}
\aliascntresetthe{rem}

\newaliascnt{ex}{theorem}
\newtheorem{ex}[ex]{Example}
\aliascntresetthe{ex}

\def\b{{\mathcal{B}}}
\def\C{{\mathbb{C}}}

\def\G{{\overline{G}}}
\def\I{{\mathfrak{I}}}

\def\N{{\mathbb{N}}}
\def\O{{\mathcal{O}}}
\def\Q{{\mathbb{Q}}}
\def\R{{\mathbb{R}}}

\def\U{{\mathcal{U}}}
\def\Z{{\mathbb{Z}}}
\def\z{{\mathcal{Z}}}

\DeclareMathOperator{\dr}{\dim_{\mathrm{Rok}}}
\DeclareMathOperator{\Aut}{\mathrm{Aut}}
\DeclareMathOperator{\id}{\mathrm{id}}
\DeclareMathOperator{\Ad}{\mathrm{Ad}}
\DeclareMathOperator{\Inf}{\mathrm{Inf}}
\DeclareMathOperator{\Hom}{\mathrm{Hom}}
\DeclareMathOperator{\scd}{\mathrm{scd}}

\hypersetup{
colorlinks,
linkcolor=blue,          
filecolor=magenta,      
urlcolor=cyan           
}

\begin{document}
\title{Values of Rokhlin Dimension for Integer Actions}
\author{Prahlad Vaidyanathan}
\email{prahlad@iiserb.ac.in}
\begin{abstract}
For an action of the integers on a C*-algebra, we show that the only possible values for the Rokhlin dimension (with commuting towers) are $0,1,2$ and $+\infty$. We study these values for certain classes of C*-algebras, and show that for any strongly outer automorphism of the Jiang--Su algebra, its Rokhlin dimension with and without commuting towers are distinct.
\end{abstract}

\maketitle
\parindent 0pt

\section{Introduction}

The notion of Rokhlin dimension was introduced by Hirshberg, Winter and Zacharias in \cite{hirshberg2015}, and may be thought of as a version of covering dimension for group actions. It is modelled on the Rokhlin property, which was introduced by Connes \cite{connes1975} in the study of automorphisms of von Neumann algebras, and is itself inspired by the Rokhlin lemma in ergodic theory. The Rokhlin property was then adapted to C*-algebras by Herman and Ocneanu \cite{herman1984}, Kishimoto \cite{kishimoto1995}, Bratteli, Evans and Kishimoto \cite{bratteli1995}, Nakamura \cite{nakamura2000}, Izumi \cite{izumi2004} and many others. Rokhlin dimension is a higher dimensional analogue of the Rokhlin property, designed so that the zero dimensional case corresponds to the Rokhlin property. \\

Broadly speaking, Rokhlin dimension is used to prove that certain regularity properties pass from the underlying algebra to the corresponding crossed product C*-algebra. These properties include finiteness of nuclear dimension, ideal separation, $\mathcal{D}$-stability, and others that are useful in the classification programme. \\

Rokhlin dimension arises in two flavours: with and without commuting towers, which we denote by $\dr^c$ and $\dr$, respectively. From the beginning, it was clear that many applications required the version with commuting towers, which is a more restrictive definition. Indeed, it has been proved by Gardella, Hirshberg and Santiago \cite{gardella2021} that for actions of compact groups, many regularity properties pass from the underlying C*-algebra to the crossed product algebra when $\dr^c$ is finite. However, several of these results are not known if we only assume that $\dr$ is finite. \\

In \cite{hirshberg2015a}, Hirshberg and Phillips recast the original commuting-towers definition for actions of a finite group $G$ in terms of a single homomorphism whose domain is the commutative C*-algebra $C(E_m(G))$, where $E_m(G)$ denotes the topological join of $m$ copies of $G$. This formulation was later expanded by Gardella \cite{gardella2017} to actions of compact groups.  We mention this because it has one important consequence from our perspective. A generalization of the Borsuk--Ulam theorem states that there is no $G$-equivariant map from $E_{m+1}(G)$ to $E_m(G)$ when $G$ is a non-trivial compact Lie group (see \cite{matousek2003} and \cite[Theorem 5.2]{hirshberg2026}). In effect, this implies that for each $m \geq 0$, the natural action $G\curvearrowright C(E_{m+1}(G))$ has Rokhlin dimension exactly $m$. Recently, Hirshberg and Phillips used equivariant K-theory in \cite{hirshberg2026} to construct actions of compact Lie groups on simple C*-algebras with arbitrarily large finite values of $\dr^c$. \\

The present article began as an attempt to determine if similar results may hold for actions of the integers, or more generally for actions of discrete, finitely generated, residually finite groups. The notion of Rokhlin dimension for actions of such groups was introduced by Szab\'{o}, Wu and Zacharias in \cite{szabo2019}. Given an action of such a group $G$, the definition of Rokhlin dimension with commuting towers may be recast in terms of a single homomorphism whose domain is the commutative C*-algebra $C(E_m(\G))$, where $\G$ denotes the profinite completion of $G$ (see \cite{self_suresh_res_finite}). Therefore, we began by trying to prove that the Rokhlin dimension of the action of $G$ on $C(E_{m+1}(\G))$ is exactly $m$. Not only is this false in general, but we managed to prove that for a large class of groups, the range of possible values that $\dr^c$ can take is surprisingly small.

\begin{mainthm}\label{thm_main}
Let $G$ be a non-trivial, discrete, finitely generated, residually finite group and suppose that $s := \scd(\G) < \infty$. Then, for any C*-algebra $A$ and any action $\alpha : G\to \Aut(A)$, \[
  \dr^c(\alpha) \in \{0,1,\ldots, s\}\cup\{+\infty\}.
\]
In particular, the only possible values for integer actions are $0,1,2,$ and $+\infty$.
\end{mainthm}
The number $\scd(\G)$ that occurs in this theorem is the strict cohomological dimension of $\G$, which we describe in \autoref{sec_group_cohomology}. \\

We believe that \autoref{thm_main} sheds light on a number of existing examples that are known in the literature. In \cite{kishimoto1995}, Kishimoto classified strongly outer automorphisms of UHF algebras by showing that any such automorphism satisfies the Rokhlin property. In \cite{nakamura2000}, Nakamura proved an analogous result for Kirchberg algebras. Similar results have since been found in other contexts by Liao \cite{liao2016} and Gardella, Hirshberg and Vaccaro \cite{gardella2022} as well. In each of these cases, an outerness hypothesis yields finite Rokhlin dimension with an explicit bound of $1$ or $2$. \autoref{thm_main} explains why such small values are to be expected, at least in the commuting-tower setting. \\

In the last part of this article, we try to put some of these examples into context by studying integer actions on Kirchberg algebras, UHF algebras and the Jiang--Su algebra $\z$. We also determine when certain specific values occur, which helps fine-tune our understanding of automorphisms of these algebras. In the process, we also show that $\dr(\gamma) \neq \dr^c(\gamma)$ for any strongly outer automorphism $\gamma$ of $\z$. Examples of this sort have been known for actions of finite groups (see \cite[Example 4.29]{gardella2021}, for instance). As far as we know, \autoref{thm_jiang_su_drc} is the first such example for integer actions.

\subsection{Notational Conventions} Throughout this article, $A$ will denote a C*-algebra, and $G$ will denote a finitely generated, discrete group (which will also usually be residually finite). We will use the letter $P$ to denote a profinite group, and we will use $M$ and $N$ to denote modules. The letters $X,Y$, and $Z$ will be reserved for topological spaces (which will usually be compact and Hausdorff). A \emph{map} between topological spaces will always mean a continuous function. If $X$ and $Y$ both carry an action of a group $G$, then a $G$-map $f:X\to Y$ will mean a continuous $G$-equivariant function.
\section{Preliminaries}

\subsection{The Profinite Completion}\label{sec_profinite_completion}

We begin with a number of preliminary notions we need going forward. In this section, we describe the profinite completion of a countable, discrete, residually finite group. This space carries a natural free action of the group, and will act as a model space for actions with finite Rokhlin dimension. \\

Recall that a profinite group is an inverse limit of a surjective inverse system of finite groups. Such groups have been well-studied over the years, and we point to \cite{ribes2010} for a thorough treatment of the subject. We will limit ourselves to those profinite groups that arise as the profinite completion of residually finite discrete groups. \\

Recall that a countable discrete group $G$ is said to be \emph{residually finite} if for each non-identity element $g \in G$, there is a (normal) subgroup $H < G$ of finite index such that $g \notin H$. For convenience, we write $\I(G)$ for the set of all normal subgroups of $G$ of finite index. This set admits a partial order by reverse inclusion (i.e. $H\leq K$ if $K\subset H$). Whenever $K\subset H$, there is a natural homomorphism $\varphi^K_H : G/K\to G/H$ given by $gK\to gH$, which makes the family $\{G/H, \varphi^K_H\}$ an inverse system of finite groups. We define the \emph{profinite completion of $G$} to be \[
  \G := \varprojlim (G/H,\varphi_{H}^K) = \{(g_HH)_{H\in \I(G)} : g_KH = g_HH \text{ whenever } K\subset H\}.
\]
Note that $\G$ is a group and a topological space as a subspace of $\prod_{H\in \I(G)} G/H$. Moreover, it is compact and totally disconnected, and admits a natural inclusion map $\iota : G\to \G$ with dense range (Note that $\iota$ is injective precisely because $G$ is residually finite). \\

For each $H \in \I(G)$, there is a natural action $\beta^H : G\curvearrowright G/H$ given by $\beta^H_t(gH) := tgH$. The maps $\varphi^K_H$ respect these actions, so we get an action $\beta : G\curvearrowright \G$ given by $\beta_t((g_HH)) := (tg_HH)$. Moreover, the inclusion map $\iota$ is $\beta$-equivariant. We write $\sigma^H : G\curvearrowright C(G/H)$ and $\sigma : G\curvearrowright C(\G)$ for the actions induced by $\beta^H$ and $\beta$, respectively. Then, it follows that \[
  (C(\G),\sigma)\cong \lim_{H\in \I(G)} (C(G/H), \sigma^H).
\]
Given a residually finite group, we will often be interested in a decreasing sequence $(K_n)\subset \I(G)$ that is cofinal in $\I(G)$. Following \cite[Definition 3.16]{szabo2019}, such a sequence will be called a \emph{dominating normal regular approximation} of $G$. Note that every finitely generated residually finite group admits a dominating normal regular approximation by \cite[Proposition 3.20]{szabo2019}. Moreover, given such a sequence $(K_n)$, we have \[
  (C(\G),\sigma) \cong \lim_{n\to \infty} (C(G/K_n), \sigma^{K_n}).
\]

\subsection{The Topological Join}\label{sec_topological_join}
Given two compact Hausdorff spaces $X$ and $Y$, the \emph{topological join} of $X$ and $Y$ is defined as \[
  X\ast Y := [0,1]\times X\times Y/\sim
\]
where $(0,x,y)\sim (0,x',y)$ and $(1,x,y) \sim (1,x,y')$ for all $x,x'\in X$ and $y,y' \in Y$. Given three compact Hausdorff spaces $X,Y,$ and $Z$ the spaces $(X\ast Y)\ast Z$ and $X\ast (Y\ast Z)$ are naturally homeomorphic. Therefore, the join operation is associative. Given a single compact Hausdorff space $X$ and an $m \in \N$, we define $X^{\ast (m)}$ to be the $m$-fold join of $X$ with itself. \\

Now, if $G$ is a group that acts on $X$ and $Y$, there is a natural action of $G$ on $X\ast Y$ given by $g\cdot [t,x,y] := [t,gx,gy]$. In particular, there is an action of $G$ on $X^{\ast (m)}$. Note that this action of $G$ on $X^{\ast (m)}$ is free if the action on $X$ is free. \\ 

When $G$ is a compact group, we write $E_m(G)$ for the space $G^{\ast (m)}$ as is customary. By our discussion, we see that $E_m(G)$ carries a free action of $G$, which we denote by $\beta^{(m)}$. When $G$ is residually finite, there is also a natural (free) action $\beta$ of $G$ on $\G$. With a slight abuse of notation, the induced action of $G$ on $E_m(\G)$ will also be denoted by $\beta^{(m)}$. \\

When $G$ is a finite group, this space also admits the structure of a CW-complex, which we now describe.

\begin{ex}\label{ex_complex_eng} Let $G$ be a finite group, fix $m \in \N$ and consider $X := E_m(G)$.\\
  
  Write a point in $X$ as $\langle t_ig_i\rangle = \langle t_0g_0, t_1g_1,\ldots, t_{m-1}g_{m-1}\rangle$, where $t_i\geq 0, \sum_{i=0}^{m-1} t_i = 1$, and $g_i \in G$. Written this way, $g_i$ is ignored if $t_i=0$. The action of $G$ on $X$ is given by $g\cdot \langle t_ig_i\rangle = \langle t_i(gg_i)\rangle$. The $j$-skeleton is the set of points with at most $(j+1)$ non-zero coordinates, \[
    X^j = \{\langle t_ig_i\rangle : |\{i : t_i > 0\}| \leq (j+1)\}.
  \]
  Observe that $X^0 \subset X^1\subset \ldots \subset X^{m-1} = X$. We now explain precisely what the characteristic maps and cells are. \\

  For $j\geq 0$, we write the standard $j$-simplex as $\Delta^j = \{(s_0,s_1,\ldots, s_j): s_r\geq 0, \sum s_r = 1\}$, and $\partial \Delta^j \cong S^{j-1}$. For $0\leq r\leq j$, let $\delta_r : \Delta^{j-1} \to \Delta^j$ insert a $0$ at position $r$. The face $\{s_r=0\}$ of $\Delta^j$ is then the image of $\delta_r$. Then, $X^0$ is the set of vertices: \[
    X^0 = \bigsqcup_{i=0}^{m-1} G\times \Delta^0.
  \]
  The point $(g,i) \in X^0$ corresponds to the vertex $v_{g,i} := \langle 1\cdot g \text{ at } i\rangle$. \\
  
A free $G$-orbit of $1$-cells can be described as follows: For each pair $\{i_0,i_1\}$ with $0\leq i_0 < i_1\leq m-1$ and each $k \in G$, we get a characteristic map $\varphi_{\{i_0,i_1\},k} : G\times \Delta^1 \to X$ given by \[
    \varphi_{\{i_0,i_1\},k}(g;1-s,s) := \langle (1-s)g \text{ at } i_0, sgk \text{ at } i_1\rangle.
  \]
  There are ${m\choose 2}|G|$ orbits and $X^1$ is a complete $m$-partite graph with $m$ parts of size $|G|$. \\
  
For $j\geq 2$, the pattern continues; every orbit of $j$-cells can be described as follows: For a fixed $g \in G$, a subset $S = \{0 \leq i_0 < i_1 < \ldots < i_j\leq (m-1)\}$ and a tuple $k = (k_1,k_2,\ldots, k_j) \in G^j$, the $j$-simplex $g\times \Delta^j$ is mapped onto the simplex spanned by $v_{g,i_0}, v_{gk_1,i_1}, \ldots, v_{gk_j,i_j}$. Observe that there are ${m\choose j+1}|G|^j$ free orbits of $j$-cells. \\

Finally, the last skeleton is $X^{m-1} = X$. Moreover, the action of $G$ on $E_m(G)$ is cellular in the sense of \cite[Proposition II.1.15]{dieck1987}. Indeed, $G$ permutes the $j$-cells freely.
\end{ex}

\begin{defn}\label{defn_milnor}
Given a finite group $G$, we define $E(G) := \varinjlim E_m(G)$ to be Milnor's infinite join construction from \cite{milnor1956}, equipped with the CW topology and the natural free action of $G$.
\end{defn}

Vertices in $E(G)$ are points of the form $v_{g,i}$ with $g \in G$ and $i\geq 0$. For each choice of indices $i_0 < i_1 < \ldots < i_j$ and elements $g_0, g_1, \ldots, g_j \in G$, the vertices $\{v_{g_0,i_0}, \ldots, v_{g_j,i_j}\}$ span a $j$-simplex which we denote by $[v_{g_0,i_0}, v_{g_1,i_1},\ldots, v_{g_j,i_j}]$. This simplex is naturally ordered by $i$, so $E(G)$ is a simplicial complex. Therefore, its cellular chain complex (from \autoref{defn_cellular_chain_complex}) is the same as its ordered simplicial chain complex (see \cite[Section 2.2]{hatcher2002}). \\

Finally, note that $E(G)$ carries the CW-topology. Therefore, every compact subset of $E(G)$ lies in $E_m(G)$ for some $m \in \N$. Since $E_m(G)$ is $(m-2)$-connected, $E(G)$ is weakly contractible. By Whitehead's theorem \cite[Theorem 4.5]{hatcher2002}, it follows that $E(G)$ is contractible. All these facts will be relevant to us later.

\subsection{Rokhlin Dimension}

In this section, we will review the definition of Rokhlin dimension for actions of residually finite groups, due to Szab\'{o}, Wu and Zacharias \cite{szabo2019}. We will then focus on the notion of Rokhlin dimension with commuting towers, and describe an equivalent condition for the finiteness of this notion from \cite{self_suresh_res_finite}. 

\begin{defn}
Given a C*-algebra $A$, let $\ell^{\infty}(\N,A)$ be the space of all bounded sequences in $A$ and $c_0(\N,A)$ be the subspace of sequences that vanish at infinity. If $A_{\infty} := \ell^{\infty}(\N,A)/c_0(\N,A)$, then $A$ embeds in $A_{\infty}$  as the set of all constant sequences so we identify $A$ with its image in $A_{\infty}$. For a C*-subalgebra $D\subset A$, we define 
\begin{equation*}
\begin{split}
A_\infty\cap D' &= \{ x\in A_{\infty} : xd=dx\text{ for all }d\in D \} \text{ and }\\
\mathrm{Ann}(D,A_\infty) &= \{ x\in A_{\infty} : xd=dx=0\text{ for all }d\in D\}. 
\end{split}
\end{equation*}
$\mathrm{Ann}(D,A_{\infty})$ is an ideal in $A_{\infty}\cap D'$, so we write
\[
F(D,A) := (A_{\infty} \cap D')/\mathrm{Ann}(D, A_{\infty})
\]
and $\kappa_{D,A} : A_{\infty}\cap D'\to F(D,A)$ for the corresponding quotient map. When $D = A$, we write $F(A)$ for $F(A,A)$ and $\kappa_A$ for $\kappa_{A,A}$. Note that $F(D,A)$ is unital if $D$ is $\sigma$-unital.
\end{defn}

Let $G$ be a discrete group and $\alpha : G\to \Aut(A)$ be an action of $G$ on a C*-algebra $A$. If $D$ is an $\alpha$-invariant subalgebra of $A$, there is a natural induced action of $G$ on $A_{\infty}$ and on $F(D,A)$ which we denote by $\alpha_{\infty}$ and $\widetilde{\alpha}_{\infty}$ respectively. Moreover, there is a $G$-equivariant $\ast$-homomorphism $(F(D,A)\otimes_{\max} D, \widetilde{\alpha}_{\infty}\otimes \alpha\lvert_D) \to (A_{\infty}, \alpha_{\infty})$ given on elementary tensors by $\kappa_{D,A}(x)\otimes a \to x\cdot a$. Under this $\ast$-homomorphism $1_{F(D,A)}\otimes a$ is mapped to $a$ for all $a\in D$, so we think of it as a way to multiply elements of $F(D,A)$ with elements of $D$ to obtain elements of $A_{\infty}$ (in a way that respects the action of $G$). \\

Two elements $a$ and $b$ in a C*-algebra are said to be orthogonal (in symbols, we write $a\perp b$) if $ab = ba = a^{\ast}b = ab^{\ast}=0$. A contractive, completely positive map $\varphi : A\to B$ between two C*-algebras is said to have order zero if $\varphi(a)\perp \varphi(b)$ whenever $a\perp b$. As is customary, we will use the abbreviation `c.c.p' for `contractive and completely positive'.

\begin{defn}\cite[Definitions 5.4 and 5.8]{szabo2019}\label{defn_dr}
Let $A$ be a C*-algebra, $G$ be a discrete, countable, residually finite group, $H$ be a subgroup of $G$ of finite index, and let $\alpha : G\to \Aut(A)$ be an action of $G$ on $A$. We say that $\alpha$ has \emph{Rokhlin dimension $d$ relative to $H$} if $d$ is the least non-negative integer such that for any separable, $\alpha$-invariant C*-subalgebra $D\subset A$, there exist $(d+1)$ equivariant c.c.p. order zero maps \[
  \varphi_0, \varphi_1, \ldots, \varphi_d : (C(G/H), \sigma^H)\to (F(D,A),\tilde{\alpha}_{\infty})
\]
such that $\sum_{\ell=0}^d \varphi_{\ell}(1_{C(G/H)}) = 1_{F(D,A)}$. We denote the Rokhlin dimension of $\alpha$ relative to $H$ by $\dr(\alpha, H)$. If no such integer exists, then we write $\dr(\alpha, H) = +\infty$. \\

We define the \emph{Rokhlin dimension of $\alpha$} as
\[
\dr(\alpha) = \sup\{ \dr(\alpha, K):  K < G \text{ of finite index}\}.
\]
\end{defn}

\begin{defn}\cite[Definition 10.2]{szabo2019}
Let $A$ be a C*-algebra, $G$ be a discrete, countable, residually finite group, $H$ be a subgroup of $G$ of finite index, and $\alpha:G\to \Aut(A)$ be an action of $G$ on $A$. We say that $\alpha$ has \emph{Rokhlin dimension $d$ with commuting towers relative to $H$} if $d$ is the least non-negative integer such that for any separable, $\alpha$-invariant C*-subalgebra $D\subset A$, there exist $(d+1)$ equivariant c.c.p. order zero maps
\[
\varphi_0, \varphi_1, \ldots, \varphi_d :(C(G/H),\sigma^H)\to (F(D,A),\widetilde{\alpha}_{\infty})
\]
with pairwise commuting ranges such that $\sum_{\ell=0}^{d}\varphi_{\ell}(1_{C(G/H)})=1_{F(D,A)}$. We denote the Rokhlin dimension of $\alpha$ with commuting towers relative to $H$ by $\dr^c(\alpha, H)$. If no such integer $d$ exists, then we write $\dr^c(\alpha, H) = +\infty$. \\

We define the \emph{Rokhlin dimension of $\alpha$ with commuting towers} as
\[
\dr^c(\alpha) = \sup\{ \dr^c(\alpha, H):  H < G\text{ of finite index}\}.
\]
\end{defn}
If $K< G$ has finite index and $H \in \I(G)$ is contained in $K$, then by \cite[Lemma 5.7]{szabo2019}, $\dr(\alpha,K)\leq \dr(\alpha,H)$ and $\dr^c(\alpha, K)\leq \dr^c(\alpha, H)$. Hence, \begin{equation*}
\begin{split}
  \dr(\alpha)  &= \sup\{\dr(\alpha,H) : H\in \I(G)\}, \text{ and } \\
  \dr^c(\alpha) &= \sup\{ \dr^c(\alpha, H): H \in \I(G)\}.
\end{split}
\end{equation*}
These are the definitions we will use going forward. \\

Given a residually finite group as above, we now rephrase these definitions in terms of the profinite completion $\G$ of $G$. Recall from \autoref{sec_profinite_completion} that there is a natural action of $G$ on $\G$, which we denoted by $\beta$. Indeed, for $m\geq 1$, there is a natural (free) action of $G$ on $E_m(\G)$ which we denoted by $\beta^{(m)}$. The induced action of $G$ on $C(E_m(\G))$ will be denoted by $\sigma^{(m)}$. When $m=1$, we simply write $\sigma = \sigma^{(1)}$.

\begin{prop}\cite[Proposition 2.9]{self_suresh_res_finite}\label{prop_dr_equivalent_definition}
Let $A$ be a C*-algebra, $G$ be a discrete, finitely generated, residually finite group and $\alpha:G\to \Aut(A)$ be an action of $G$ on $A$. Then $\dr^c(\alpha) \leq d$ if and only if, for each separable, $\alpha$-invariant C*-subalgebra $D\subset A$, there exist $(d+1)$ equivariant, order zero, c.c.p. maps
\[
\varphi_0, \varphi_1, \ldots, \varphi_d : (C(\G), \sigma) \to (F(D,A),\widetilde{\alpha}_{\infty})
\]
with pairwise commuting ranges such that $\sum_{\ell=0}^d \varphi_{\ell}(1_{C(\G)}) = 1_{F(D,A)}$.
\end{prop}

\begin{theorem}\cite[Theorem 2.11]{self_suresh_res_finite}\label{thm_drc_equivalent_definition}
Let $A$ be a C*-algebra, $G$ be a finitely generated, residually finite group and $\alpha:G\to \Aut(A)$ be an action of $G$ on $A$. Then, $\dr^c(\alpha) \leq d$ if and only if, for each separable, $\alpha$-invariant C*-subalgebra $D\subset A$, there is a unital, $G$-equivariant $\ast$-homomorphism
\[
\varphi : (C(E_{d+1}(\G)), \sigma^{(d+1)}) \to (F(D,A), \widetilde{\alpha}_{\infty}).
\]
\end{theorem}
Indeed, an analogous statement holds for subgroups as well. If $H \in \I(G)$, then $\dr^c(\alpha,H)\leq d$ if and only if, for each separable, $\alpha$-invariant C*-subalgebra $D\subset A$, there is a unital $G$-equivariant $\ast$-homomorphism \[
  \varphi : (C(E_{d+1}(G/H)), (\sigma^H)^{(d+1)}) \to (F(D,A),\tilde{\alpha}_{\infty}).
\]
The proof of this is similar to the proof of \cite[Theorem 2.11]{self_suresh_res_finite}, so we omit it. We will, however, use this fact tacitly in the future.

\subsection{Obstruction Theory}\label{sec_obstruction_theory}

We now turn our attention to an idea that is of crucial importance to us. Obstruction theory is a classical subject that concerns itself with the extension problem: Given a relative CW-complex $(X,A)$ and a (continuous) map $f_A : A\to Y$, it is used to determine if there is a map $f:X\to Y$ that extends $f_A$. In this paper, we will be concerned with the equivariant version of this problem: if we are given a finite group $G$ and two $G$-spaces $X$ and $Y$, we would like to know if we can construct a continuous $G$-map $f:X\to Y$. \\

For this section, we follow the exposition of \cite{dieck1987}. To begin to explain the theory, we need to define equivariant cohomology groups. Note that throughout this section, $G$ will denote a finite group and $\Z[G]$ will denote the associated group ring. By a $G$-module, we will mean a module over $\Z[G]$. For two $G$-spaces $X$ and $Y$, we write $[X,Y]_G$ for the $G$-homotopy classes of maps between them. When $(X,A)$ is a relative $G$-complex, we will write $[X,Y]_G$ for the homotopy classes $\mathrm{rel}A$ of extensions of a fixed $G$-map $f_A : A\to Y$. In many of our applications, we will take $A$ to be the empty set, so this abuse of notation should not be confusing. Finally, we write $[X,Y]$ for the homotopy classes of continuous maps from $X$ to $Y$.

\begin{defn}\label{defn_relative_g_complex}
Let $G$ be a finite group, and let $(X,A)$ be a relative $G$-complex with a free $G$-action on $X\setminus A$. In other words, we have an equivariant CW-decomposition $(X^n)$ of $(X,A)$ where ~\begin{enumerate}
    \item $X^{-1} := A \subset X^0$ and $X = \bigcup X^n$.
    \item For each $n\geq 0$, $X^n$ is obtained from $X^{n-1}$ by attaching a disjoint union of free $G$-cells of the form $G\times D^n$.
    \item $X$ carries the colimit topology with respect to $(X^n)$. i.e. $B\subset X$ is closed if and only if $B\cap X^n$ is closed for each $n\in \N_0$.
  \end{enumerate}
\end{defn}
\begin{defn}\label{defn_cellular_chain_complex}
Now consider the (singular) cellular chain complex $C_{\ast}(X,A)$ given by \[
    \ldots \to H_k(X^k,X^{k-1})\xrightarrow{d} H_{k-1}(X^{k-1},X^{k-2})\to \ldots
  \]
  Each $H_k(X^k,X^{k-1})$ is a free $G$-module, and the boundary map $d : H_k(X^k,X^{k-1})\to H_{k-1}(X^{k-1},X^{k-2})$ is also $G$-equivariant by naturality of the connecting homomorphism of the triple $(X^k,X^{k-1},X^{k-2})$. Therefore, $C_{\ast}(X,A)$ becomes a chain complex of $G$-modules. \\

  Given a $G$-module $M$, define the cochain complex $C_G^{\ast}(X,A,M) := \Hom_{\Z[G]}(C_{\ast}(X,A),M)$. The equivariant cohomology groups are defined as the corresponding homology groups, \[
    H_G^m(X,A,M) := H_m(C_G^{\ast}(X,A,M)).
  \]
  If $A = \emptyset$, we simply write $H^m_G(X,M)$ for this group.  
\end{defn}

Now fix $n\geq 1$ and let $Y$ be a path-connected, $n$-simple space equipped with a (not necessarily free) $G$-action. (Recall that $Y$ is $n$-simple if $\pi_1(Y,y)$ acts trivially on $\pi_n(Y,y)$). Then, the forgetful map $\pi_n(Y) \to [S^n,Y]$ is bijective, and hence $G$ acts on $\pi_n(Y)$. We will then be interested in the cohomology groups $H_G^m(X,A,\pi_n(Y))$. \\

The first important theorem for our context is the following result. Given a $G$-map $h : X^n\to Y$, we wish to extend $h\lvert_{X^{n-1}}$ to a $G$-map $\tilde{h} : X^{n+1}\to Y$ (after potentially changing $h$ on $n$-cells). The obstruction to such an extension, called the \emph{obstruction class}, is an element $c_G^{n+1}(h) \in H_G^{n+1}(X,A,\pi_n(Y))$ and has the property that $c_G^{n+1}(h) = 0$ if and only if such an extension $\tilde{h}$ exists. \\

Moreover, the class $c_G^{n+1}(h)$ has an explicit description as a cellular cochain (which we denote by $\mathfrak{c}_G^{n+1}(h)$): Let $\phi = (\phi_j) : \bigsqcup_j G\times (D^{n+1},S^n)\to (X^{n+1},X^n)$ be a characteristic map and $e_j \in C_{n+1}(X,A)$ be the basis element corresponding to $\phi_j$, then $\mathfrak{c}_G^{n+1}(h)(e_j) \in [S^n,Y]$ is the homotopy class represented by $h\phi_j\lvert_{\{1\}\times S^n}$.

\begin{theorem}\cite[Theorem II.3.10]{dieck1987}\label{thm_obstruction_exact_sequence}
For each $n\geq 1$, there is an exact obstruction sequence \begin{equation*}
  [X^{n+1},Y]_G \to \mathrm{Im}([X^n,Y]_G\to [X^{n-1},Y]_G)\xrightarrow{c^{n+1}_G} H_G^{n+1}(X,A,\pi_n(Y)).
\end{equation*}
\end{theorem}

In the context we will be interested in, $Y$ will not only be an $n$-simple space, but also $(n-1)$-connected. In that case, the extension problem simplifies even further.

\begin{prop}\label{prop_extension_uniqueness}\cite[Proposition II.3.15]{dieck1987}
Suppose $Y$ is $n$-simple and $(n-1)$-connected, and suppose that $(X,A)$ is a relative $G$-complex with a free action on $X\setminus A$. Then, the following hold: ~\begin{enumerate}
  \item Every $G$-map $f:A\to Y$ is extendable to a $G$-map $\tilde{f} : X^n\to Y$. Moreover, any two extensions are homotopic $\mathrm{rel}A$ on $X^{n-1}$.
  \item Let $k:f_0\simeq f_1 : A\to Y$ be a homotopy and let $F_0, F_1 : X^n\to Y$ be extensions of $f_0$ and $f_1$ respectively. Then, there is a homotopy $K$ connecting $F_0\lvert_{X^{n-1}}\simeq F_1\lvert_{X^{n-1}}$ and extending $k$.
\end{enumerate}
\end{prop}

\begin{rem}\label{rem_primary_obstruction} Let $G$ be a finite group acting freely on a relative $G$-complex $(X,A)$ and assume $Y$ is $n$-simple and $(n-1)$-connected. Now suppose $f_A:A\to Y$ is fixed, then by \autoref{prop_extension_uniqueness}, there is an extension $f : X^n\to Y$ and that determines a unique element in \[
    \mathrm{Im}([X^n,Y]_G \to [X^{n-1},Y]_G).
  \]
  Hence, we have a well-defined element $c^{n+1}_G(f) \in H_G^{n+1}(X,A,\pi_n(Y))$. By \autoref{prop_extension_uniqueness}, we obtain a well-defined map \[
    \gamma_{n+1} : [A,Y]_G \to H_G^{n+1}(X,A,\pi_n(Y)) \text{ given by } \gamma_{n+1}([f_A]) := c^{n+1}_G(f).
  \]
  Moreover, $\gamma_{n+1}([f_A]) = 0$ if and only if $f_A$ has an extension to $X^{n+1}$, and therefore, $\gamma_{n+1}([f_A])$ is called the \emph{primary obstruction to extending $f_A$}. \\

  Now taking $A = \emptyset$, there is a unique map $f_A:\emptyset \to Y$, so we get an element \[\gamma_{n+1}(X,Y) \in H_G^{n+1}(X,\pi_n(Y)).\] 
  By \autoref{prop_extension_uniqueness}, any two $G$-maps from $X^n$ to $Y$ agree up to homotopy on $X^{n-1}$, so $\gamma_{n+1}(X,Y) = 0$ if and only if there is a $G$-map from $X^{n+1}$ to $Y$. 
\end{rem}

Finally, we need the fact that the obstruction exact sequence from \autoref{thm_obstruction_exact_sequence} is natural, and hence so is the construction of the primary obstruction described above. Moreover, it is also natural across a homomorphism between the groups as well. We record this fact below. The proof of this result follows from the explicit description of the obstruction class given before \autoref{thm_obstruction_exact_sequence} (also see the comments on \cite[Pg. 116]{dieck1987}). 

\begin{prop}\label{prop_natural_primary_obstruction}
Let $\rho :G'\to G$ be a homomorphism between finite groups, $X$ be a free $G'$-complex, $Z$ be a free $G$-complex, and $\theta : X\to Z$ be a cellular $\rho$-equivariant map. Fix $n \in \N$, and let $Y$ be an $n$-simple and $(n-1)$-connected $G$-space (treated as a $G'$-space via $\rho$). If $\gamma_{n+1}(X,Y)$ and $\gamma_{n+1}(Z,Y)$ are the primary obstructions described above, then \[
  \gamma_{n+1}(X,Y) = \theta^{\ast}\gamma_{n+1}(Z,Y),
\]
where $\theta^{\ast} : H_G^{n+1}(Z,\pi_n(Y)) \to H_{G'}^{n+1}(X,\rho^{\ast}\pi_n(Y))$ is the natural map induced by $\theta$.
\end{prop}

\section{Main Results}\label{sec_main}
\subsection{Cohomology of $\G$}\label{sec_group_cohomology}
We now describe the continuous cohomology of the group $\G$, the profinite completion of a residually finite group $G$. Indeed, the main goal of this section is to prove an inductive limit decomposition of the cohomology groups of $\G$ in terms of the equivariant cohomology groups of $E(G/K)$, as $K$ varies over $\I(G)$ (see \autoref{cor_inductive_limit_cohomology_space}). The reader is referred to \cite[Chapter 6]{ribes2010} for some of the details, along with \cite{brown1982} for the usual cohomology of finite groups.\\

For the moment, let $G$ be a finite group and consider $\Z$ as a trivial $G$-module. Let $F_{\ast}(G)$ be the standard free resolution of $\Z$ over $\Z[G]$ (see \cite[Section I.5]{brown1982}). Explicitly, $F_n(G)$ is the free abelian group on $G^{n+1}$ with the diagonal action of $G$. The boundary map $\partial: F_n(G) \to F_{n-1}(G)$ is given by $\partial = \sum_{i=0}^n (-1)^id_i$ where $d_i(g_0,\ldots, g_n) = (g_0,\ldots, \hat{g_i},\ldots, g_n)$. We also have an augmentation map $\varepsilon : F_0(G)\to \Z$ given by $\varepsilon(g) = 1$ for all $g \in G$. Given a $G$-module $M$, define $F^n(G,M) := \Hom_G(F_n(G),M)$. The group cohomology of $G$ with coefficients in $M$, denoted $H^{\ast}(G,M)$, is the homology of this cochain complex $(F^n(G,M))$. \\

There is another natural way to define these cohomology groups. Since we need to compare these two definitions, we give an explicit description. Let $E(G)$ denote Milnor's infinite join construction from \autoref{defn_milnor}. Then, $E(G)$ is a contractible, free $G$-complex and its cellular chain complex may be identified with the ordered simplicial chain complex. \\

The latter is defined as follows: Let $C_j(E(G))$ be the free abelian group on the $j$-simplices of $E(G)$ described after \autoref{defn_milnor}, and define the boundary map $\partial : C_j(E(G))\to C_{j-1}(E(G))$ by $\partial = \sum_{s=0}^j (-1)^s d_s$ where $d_s([v_{g_0,i_0}, \ldots, v_{g_j,i_j}]) = [v_{g_0,i_0}, \ldots, \widehat{v_{g_s,i_s}}, \ldots v_{g_j,i_j}]$. Note that $G$ acts on $C_j(E(G))$ by acting on vertices by $g\cdot v_{g_s,i_s} := v_{gg_s,i_s}$. This action permutes the simplices freely, so each $C_j(E(G))$ is a free $G$-module. We define the augmentation map $\varepsilon : C_0(E(G))\to \Z$ by $\varepsilon(v_{g,i}) := 1$. This makes $C_{\ast}(E(G))\to \Z$ a free resolution of $\Z$ over $\Z[G]$. \\

Given another $G$-module $M$, we define $C_G^n(E(G),M) := \Hom_G(C_n(E(G)), M)$. The sequence $(C_G^n(E(G),M))$ is then a cochain complex, whose homology groups are the equivariant cohomology groups $H^{\ast}_G(E(G),M)$ defined in \autoref{sec_obstruction_theory}. The next lemma explicitly identifies these two constructions (see also \cite[Section I.4]{brown1982}, where these are identified abstractly).

\begin{lemma}\label{lem_cohomology_identification}
The map $\delta_G : C_j(E(G)) \to F_j(G)$ defined on basis elements by \[
  [v_{g_0,i_0}, \ldots, v_{g_j,i_j}] \mapsto (g_0,\ldots, g_j)
\]
defines an augmentation-preserving chain map of $G$-modules. Hence, the map $\delta_G^{\ast} : F^{\ast}(G,M) \to C_G^{\ast}(E(G),M)$ given by $\sigma \mapsto \sigma \circ \delta_G$ induces an isomorphism of cohomology groups, \[
  H^m(G,M) \cong H^m_G(E(G),M).
\]
\end{lemma}
\begin{proof}
It is easy to check that $\delta_G$ is a $G$-linear chain map that preserves the augmentation map. Since $C_{\ast}(E(G))$ and $F_{\ast}(G)$ are both projective resolutions of $\Z$ over $\Z[G]$, $\delta_G$ is a homotopy equivalence by \cite[Theorem I.7.5]{brown1982}. Hence, $\delta_G^{\ast}$ induces an isomorphism on cohomology.
\end{proof}

For our purposes, we will also need the continuous cohomology (with discrete coefficients) of a profinite group. Here, we replace the cochain complex $F^{\ast}(G,M)$ by its continuous counterpart. \\

Let $P$ be a profinite group and $M$ be a discrete $P$-module. For $n \geq 0$, we write $F^n(P,M)$ for the set of all \emph{continuous} functions $f : P^{n+1} \to M$ satisfying $f(gg_0,gg_1,\ldots, gg_n) = gf(g_0,g_1,\ldots, g_n)$ for all $g,g_i \in P$. There is a boundary map $\partial^{n+1} : F^n(P,M)\to F^{n+1}(P,M)$ which yields a cochain complex $F^{\ast}(P,M)$. Then the $m^{th}$ cohomology groups of $P$ with coefficients in $M$, denoted $H^m(P,M)$, are defined as the homology groups of this cochain complex. \\

In what follows, let $G$ denote a discrete, finitely generated, residually finite group and let $(K_n)$ be a fixed dominating normal regular approximation. Then, $\G$ is the inverse limit of the sequence $(G/K_n)$ and we will think of $G$ as a subgroup of $\G$ via the natural inclusion map $\iota : G\to \G$. Now, we wish to compute the cohomology groups $H^m(\G,M)$ in terms of the cohomology groups of $G/K_n$. \\

Firstly, if $K \in \I(G)$ is fixed, then the set $\{H\in \I(G) : H\subset K\}$ is a cofinal set in $\I(K)$, so $\overline{K}$ (the profinite completion of $K$) may be thought of as a subgroup of $\G$. Also, $\overline{K}$ is the closure of $K$ in $\G$, it is normal in $\G$ and $\G/\overline{K}\cong G/K$. This allows us to make the following definition.

\begin{defn}\cite[Section 6.5]{ribes2010}\label{defn_inflation_map}
Let $K \in \I(G)$ and let $M$ be a discrete $\G$-module. Then define \[
  M^K := \{m\in M : km = m \text{ for all } k\in K\}.
\]
Then, $M^{\overline{K}} = M^K$ and hence $M^K$ is a module over $\G/\overline{K}\cong G/K$. Now the projection map $\pi_K : \G\to G/K$ and the inclusion map $M^K\to M$ induce a homomorphism \[
  \mathrm{Inf}_{\G}^{G/K} : H^m(G/K,M^K)\to H^m(\G,M).
\]
This is called the \emph{inflation map}.
\end{defn}


Moreover, if $K,H \in \I(G)$ are such that $K\subset H$, then there is an inflation map \[
  \Inf^{G/H}_{G/K} : H^m(G/H,M^H)\to H^m(G/K,M^K)
\]
with the property that $\Inf_{\G}^{G/K}\Inf^{G/H}_{G/K} = \Inf_{\G}^{G/H}$.\\

Therefore, for a fixed $m \in \N$, the sequence $(H^m(G/K_n,M^{K_n}), \mathrm{Inf}_{G/K_n}^{G/K_{n+1}})$ is an inductive sequence of abelian groups. The following result is now the  expected conclusion of this discussion.

\begin{prop}\cite[Corollary 6.5.6]{ribes2010}\label{prop_inductive_limit_cohomology_module}
For each $m \geq 0$ and each discrete $\G$-module $M$, we have 
\[
  H^m(\G,M) = \lim_{n\to \infty} (H^m(G/K_n,M^{K_n}), \Inf_{G/K_n}^{G/K_{n+1}}).
\]
\end{prop}

We now wish to understand these inflation maps from the point of view of \autoref{lem_cohomology_identification}. To that end, we work with the following setup: Let $\varphi : G'\to G$ be a homomorphism of finite groups, let $N$ be a $G$-module and $N'$ be a $G'$-module, and let $j: N\to N'$ be a $\varphi$-equivariant homomorphism  of modules. \\

The map $\varphi$ then induces a map of standard resolutions $F(\varphi) : F_n(G')\to F_n(G)$ given by \[
  (g_0,g_1,\ldots, g_n)\mapsto (\varphi(g_0),\varphi(g_1),\ldots, \varphi(g_n)).
\]
It also induces a map $\tilde{\varphi} : E(G')\to E(G)$ given by $\langle t_ig_i\rangle \mapsto \langle t_i\varphi(g_i)\rangle$, which in turn induces a chain map $\tilde{\varphi}_{\ast} : C_{\ast}(E(G'))\to C_{\ast}(E(G))$ given on basis elements by \[
  [v_{g_0,i_0},\ldots, v_{g_j,i_j}]\mapsto [v_{\varphi(g_0),i_0},\ldots, v_{\varphi(g_j),i_j}].
\]
It is now obvious that \[
  \delta_G\circ \tilde{\varphi}_{\ast} = F(\varphi)\circ \delta_{G'}.
\]
Now define two maps $\eta : F^{\ast}(G,N)\to F^{\ast}(G',N')$ and $\lambda : C^{\ast}_G(E(G),N)\to C^{\ast}_{G'}(E(G'),N')$ by \[
  \eta(\sigma) := j\circ \sigma\circ F(\varphi), \text{ and } \lambda(f) := j\circ f \circ \tilde{\varphi}_{\ast}.
\]
It is clear that these are both maps of cochain complexes and that \[
  \lambda\circ \delta_G^{\ast} = \delta_{G'}^{\ast}\circ \eta.
\]
We therefore obtain the following result which will be a useful computational tool for us.

\begin{prop}\label{prop_cohomology_isomorphism}
The following diagram commutes: \[
  \xymatrix{
    H^m(G,N)\ar[rr]_{\delta_G^{\ast}}^{\cong}\ar[d]_{\eta_{\ast}} && H^m_G(E(G),N)\ar[d]^{\lambda_{\ast}} \\
    H^m(G',N')\ar[rr]_{\delta_{G'}^{\ast}}^{\cong} && H^m_{G'}(E(G'),N').
  }
\]
Moreover, if $N' = \varphi^{\ast}N$ and $j:N\to N'$ is the identity map, then $\lambda_{\ast} = \tilde{\varphi}^{\ast}$.
\end{prop}

We are now in a position to understand the inflation map under the identification of \autoref{lem_cohomology_identification}. Fix $K,H \in \I(G)$ with $K\subset H$, and a discrete $\G$-module $M$. As described in \cite[Section 6.5]{ribes2010}, the map $\Inf_{G/H}^{G/K}$ may be described at the level of homogeneous cochains as follows: Let $\varphi^K_H : G/K\to G/H$ be the natural homomorphism. Given $\sigma \in F^{n}(G/H,M^H)$, we define $\eta^K_H(\sigma)\in F^{n}(G/K,M^K)$ by 
\[
  (x_0,x_1,\ldots, x_n) \mapsto \sigma(\varphi^K_H(x_0),\varphi^K_H(x_1),\ldots, \varphi^K_H(x_{n})).
\]
This defines a map of chain complexes, \[
  \eta^K_H : F^{n}(G/H,M^H)\to F^{n}(G/K,M^K),
\]
which in turn induces the inflation map $\Inf_{G/H}^{G/K}$ at the level of homology.  \\

Hence, if we take $G' = G/K, G = G/H, \varphi = \varphi^K_H$ and $j : M^H\hookrightarrow M^K$ to be the inclusion map in \autoref{prop_cohomology_isomorphism}, then the map $\eta_{\ast}$ appearing on the left is $\Inf_{G/H}^{G/K}$. We write $(\lambda^K_H)_{\ast}$ for the map $H^m_{G/H}(E(G/H),M^H)\to H^m_{G/K}(E(G/K),M^K)$  that appears on the right. Applying this to the sequence of subgroups $(K_n)$ as above, we get the following consequence of \autoref{prop_inductive_limit_cohomology_module}.

\begin{cor}\label{cor_inductive_limit_cohomology_space}
For each $m \geq 0$ and every discrete $\G$-module $M$, we have \[
  H^m(\G,M) \cong \lim_{n\to\infty} (H^m_{G/K_n}(E(G/K_n),M^{K_n}), (\lambda^{K_{n+1}}_{K_n})_{\ast}).
\]  
In particular, if $K_0$ acts trivially on $M$, then $M^{K_n} = M$ for all $n$, and therefore \[
  H^m(\G,M) \cong \lim_{n\to\infty} (H^m_{G/K_n}(E(G/K_n),M), (\tilde{\varphi}^{K_{n+1}}_{K_n})^{\ast}),
\]
where $\tilde{\varphi}^K_H : E(G/K)\to E(G/H)$ is the natural map induced by the map $\varphi^K_H : G/K\to G/H$.
\end{cor}

Before we go further, we need one more important concept, that of (strict) cohomological dimension. For discrete groups, this is described in \cite[Chapter VIII]{brown1982}. In our case, we will need this notion for profinite groups, so we follow \cite[Chapter 7]{ribes2010}. \\

In what follows, given an abelian group $T$ and a prime number $p$, we write $T_p$ for the $p$-primary component of $T$, i.e. those elements in $T$ whose order is a power of $p$.

\begin{defn} Let $P$ be a fixed profinite group, and $p$ a prime. The \emph{strict cohomological $p$-dimension}, denoted $\scd_p(P)$, is the smallest non-negative $n$ such that $H^m(P,M)_p = 0$ for all $m > n$, and for all discrete $P$-modules $M$. We write $\scd_p(P) = +\infty$ if no such number exists. The \emph{strict cohomological dimension} of $P$ is defined as $\scd(P) := \sup_p \scd_p(P)$.
\end{defn}

\begin{rem}\label{prop_cohomology_zero}
The strict cohomological dimension is meant to measure the stage at which the cohomology groups of $P$ vanish. Indeed, if $P$ is a profinite group and $M$ is a discrete $P$-module, then $H^m(P,M)$ is a torsion abelian group for $m\geq 1$ (see \cite[Corollary 6.7.4]{ribes2010}). Every torsion abelian group is a direct sum of its $p$-primary components, so if $m > \scd(P)$, then $H^m(P,M) = 0$.
\end{rem}

\begin{lemma}\label{lem_scd_atleast_2}
For any non-trivial profinite group $P$, $\scd(P)\geq 2$.
\end{lemma}
\begin{proof}
Fix a prime $p$ dividing the (supernatural) order of $P$, so that $P$ admits a non-trivial $p$-Sylow subgroup $K$ (see \cite[Corollary 2.3.6]{ribes2010}). Fix a non-trivial $x \in K$ and consider the subgroup $C := \overline{\langle x\rangle}$. This is a closed procyclic, pro-$p$ group. Moreover, by \cite[Theorem 7.3.1]{ribes2010}, we have \[
  \scd(C)\leq \scd(P).
\]
It therefore suffices to prove that $\scd(C)\geq 2$. \\

By \cite[Proposition 2.7.1]{ribes2010}, we have two possibilities: Either $C$ is finite, or $C = \Z_p = \varprojlim \Z/p^n\Z$, the group of $p$-adic integers. \\

If $C$ is finite, then it is well-known (see, for instance, \cite[Section III.1]{brown1982}) that $H^{2k}(C,\Z)$ is non-zero for all $k \in \N$. Hence, $\scd(C) = +\infty$. \\

If $C \cong \Z_p$, then consider the short exact sequence $0 \to \Z\to \Q \to \Q/\Z\to 0$ with the trivial action. Since $\Q$ is divisible and torsion-free, we know from \cite[Corollary 6.7.5]{ribes2010} that $H^n(C,\Q) = 0$ for all $n\geq 1$. By the long exact sequence in cohomology, we have \[
  H^2(C,\Z)\cong H^1(C,\Q/\Z) \cong \Hom(C,\Q/\Z).
\]
Note that $\Hom(C,\Q/\Z)$ denotes the collection of continuous homomorphisms from $C$ to $\Q/\Z$, and its identification with $H^1(C,\Q/\Z)$ follows from \cite[Section 6.8]{ribes2010} and the fact that the action of $C$ on $\Q/\Z$ is trivial. \\

Now, we have a non-trivial map $\varphi : C\to \Q/\Z$ given by the composition \[
  C \to \Z/p\Z\to \Q/\Z,
\]
where the first map is the quotient map and the second is the map $\overline{1} \mapsto \overline{1/p}$. The map $\varphi$ is continuous because it is the composition of two continuous maps. Hence, as discussed above, \[
  H^2(C,\Z) \neq 0,
\]
and therefore $\scd(C)\geq 2$. This completes the proof.
\end{proof}

For our purposes, the profinite group in question will be the profinite completion $\G$ of a countable, discrete, residually finite group $G$. In this context, a few values of the strict cohomological dimension are known.

\begin{ex}
~\begin{enumerate}
  \item If $G = \Z$, then $\scd(\G) = 2$ by \cite[Example 7.1.3]{ribes2010}. This will be the most important example for us in \autoref{sec_integer_actions}.
  \item Let $G = \mathbb{F}_n$ be a free group on finitely many generators. Then $\G$ is a free profinite group by \cite[Proposition 3.3.6]{ribes2010}. By \cite[Corollary 7.5.2, Proposition 7.1.2]{ribes2010}, $\scd(\G) \leq 2$. By \autoref{lem_scd_atleast_2}, we conclude that $\scd(\G) = 2$.
  \item If $G = \Z^n$, then $\scd(\G) \leq n+1$ by \cite[Proposition 7.4.2]{ribes2010}.
  \item If $G$ is a non-trivial finite group, then $\scd(\G) = \scd(G) = +\infty$ by \cite[Corollary VIII.2.5]{brown1982}.
\end{enumerate}
\end{ex}

\subsection{Proof of \autoref{thm_main}}
We now have all the pieces in place to prove the main theorems of this article. The following theorem is best understood in the context of Dold's generalization of the Borsuk--Ulam theorem (see \cite{matousek2003} and \cite[Theorem 5.2]{hirshberg2026}). That theorem states that for a non-trivial finite group $G$, there is no $G$-map $f : E_{n+2}(G)\to E_{n+1}(G)$. For $n\geq 2$, the non-existence of such a map may be proved using equivariant obstruction theory. Indeed, such a map does not exist because the primary obstruction element \[
  \gamma_{n+1}(E_{n+2}(G),E_{n+1}(G)) \in H_G^{n+1}(E_{n+2}(G),\pi_{n}(E_{n+1}(G)))
\]
does not vanish. \\

In the following theorem, we fix a non-trivial residually finite group $G$ and two subgroups $K,H \in \I(G)$ such that $K\subset H$. We then ask if there is a $G$-map $f : E_{n+2}(G/K)\to E_{n+1}(G/H)$. Once again, the existence of such a map is determined by the vanishing of an obstruction element \[
  \gamma_{n+1}(E_{n+2}(G/K), E_{n+1}(G/H)) \in H_{G/K}^{n+1}(E_{n+2}(G/K), \pi_n(E_{n+1}(G/H))).
\]
In this case, however, this obstruction element can vanish for any given $H \in \I(G)$, provided we choose the subgroup $K$ carefully. This leads to the following result.

\begin{theorem}\label{thm_borsuk_ulam_not}
Let $G$ be a non-trivial, discrete, finitely generated, residually finite group and let $H \in \I(G)$. If $n \geq \scd(\G)$, then there is a $K \in \I(G)$ contained in $H$ which admits a $G$-map \[
  f : E_{n+2}(G/K)\to E_{n+1}(G/H).
\]
\end{theorem}
\begin{proof}
Let $Y := E_{n+1}(G/H)$. Since $n\geq 2$ (by \autoref{lem_scd_atleast_2}) and $Y$ is $(n-1)$-connected, $Y$ is simply connected and hence $n$-simple. Moreover, $M := \pi_n(Y) = [S^n,Y]$ as mentioned before, and $\G$ acts on $M$ via the quotient map $\G\to G/H$. Hence, $M$ is a discrete $\G$-module. Since $n\geq \scd(\G)$,  \[
  H^{n+1}(\G,M) = 0.
 \]
by \autoref{prop_cohomology_zero}. \\

Fix a dominating normal regular approximation $(K_m)$ of $G$ contained in $H$, with $H = K_0$. Write $X_{m} := E_{n+2}(G/K_m)$, then $X_m$ is a free $G/K_m$-complex of dimension $(n+1)$ by \autoref{ex_complex_eng}. Therefore, by \autoref{rem_primary_obstruction}, there is a $G$-map $f_m : X_m\to Y$ if and only if the primary obstruction \[\gamma_{n+1}(X_m,Y) \in H_{G/K_m}^{n+1}(X_m,M)\] vanishes. \\

Consider the natural inclusion map $\iota_{K_m} : X_m\to E(G/K_m)$, which is $G/K_m$-equivariant and cellular. This induces a map \[\iota_{K_m}^{\ast} : H^{n+1}_{G/K_m}(E(G/K_m),M)\to H^{n+1}_{G/K_m}(X_m,M).\]
Now, $E(G/K_m)$ is a free $G/K_m$-complex by the discussion following \autoref{defn_milnor}. Hence, there is a primary obstruction element $\gamma_{n+1}(E(G/K_m),Y) \in H^{n+1}_{G/K_m}(E(G/K_m),M)$ with the property that $\gamma_{n+1}(E(G/K_m),Y) = 0$ if and only if there is a $G$-map from the $(n+1)$-skeleton of $E(G/K_m)$ to $Y$. Moreover, by \autoref{prop_natural_primary_obstruction}, we have \[
  \gamma_{n+1}(X_{m},Y) = \iota_{K_m}^{\ast}\gamma_{n+1}(E(G/K_m),Y).
\]
 For each $m \geq \ell \geq 0$, write $\varphi^{K_{m}}_{K_{\ell}} : G/K_{m} \to G/K_{\ell}$ for the corresponding quotient map, and write \[
  \tilde{\varphi}^{K_m}_{K_{\ell}} : E(G/K_m)\to E(G/K_{\ell})
 \]
 for the induced map on spaces. Then by \autoref{prop_natural_primary_obstruction}, \[
  \gamma_{n+1}(E(G/K_{m}),Y) = (\tilde{\varphi}^{K_{m}}_{K_{\ell}})^{\ast}\gamma_{n+1}(E(G/K_{\ell}),Y) = (\tilde{\varphi}^{K_{m}}_{H})^{\ast}\gamma_{n+1}(E(G/H),Y).
\]
As mentioned above, $H^{n+1}(\G,M) = 0$. Moreover, $K_0 = H$ acts trivially on $M$, so by \autoref{cor_inductive_limit_cohomology_space},
\[
  \lim_{m\to\infty} (H^{n+1}_{G/K_m}(E(G/K_m),M),(\tilde{\varphi}^{K_{m+1}}_{K_m})^{\ast}) = 0.
\]
Hence, there is $m \in \N$ such that \[
  \gamma_{n+1}(E(G/K_m),Y) = (\tilde{\varphi}^{K_{m}}_{H})^{\ast}\gamma_{n+1}(E(G/H),Y) = 0.
\]
It follows that \[\gamma_{n+1}(X_{m},Y) = \iota_{K_m}^{\ast}\gamma_{n+1}(E(G/K_m),Y) = 0.\] Therefore, there is a $G$-map $f_m : X_{m}\to Y$, so the subgroup $K = K_m$ has the desired property.
\end{proof}

The proof of \autoref{thm_main} is now immediate.

\begin{proof}[Proof of \autoref{thm_main}]
Let $G$ be a non-trivial, discrete, finitely generated, residually finite group such that $s := \scd(\G) < \infty$, and let $\alpha : G\to \Aut(A)$ be an action with $d := \dr^c(\alpha) < \infty$. \\

Suppose that $d\geq s+1$, then fix a separable, $\alpha$-invariant subalgebra $D\subset A$. By \autoref{thm_drc_equivalent_definition}, there is a unital $G$-equivariant $\ast$-homomorphism \[
  \varphi : (C(E_{d+1}(\G)),\sigma^{(d+1)})\to (F(D,A), \tilde{\alpha}_{\infty}).
\]
Fix a subgroup $H\in \I(G)$. We wish to prove that $\dr^c(\alpha,H)\leq s$. Since $s = \scd(\G)$, we may iteratively apply \autoref{thm_borsuk_ulam_not} to obtain $K_i \in \I(G), i\in \{1,2,\ldots, (d-s)\}$ with $K_{d-s}\subseteq \ldots \subseteq K_1\subseteq H$, together with $G$-maps \[
  E_{d+1}(G/K_{d-s})\to \ldots \to E_{s+2}(G/K_1)\to E_{s+1}(G/H).
\]
We write $f : E_{d+1}(G/K_{d-s})\to E_{s+1}(G/H)$ for the composition of these maps. Composing with the natural map $E_{d+1}(\G)\to E_{d+1}(G/K_{d-s})$ and dualizing, we get a $G$-equivariant unital $\ast$-homomorphism \[
  \psi : C(E_{s+1}(G/H))\to C(E_{d+1}(\G)).
\]
Composing with $\varphi$ gives a $G$-equivariant unital $\ast$-homomorphism $\varphi\circ \psi : C(E_{s+1}(G/H))\to F(D,A)$. By the note after \autoref{thm_drc_equivalent_definition} (and since $D$ was arbitrary), this implies that $\dr^c(\alpha,H) \leq s$. This holds for each $H \in \I(G)$, therefore $\dr^c(\alpha)\leq s$.
\end{proof}

\section{Values for Integer Actions}\label{sec_integer_actions}

We now turn our attention to some examples of single automorphisms (identified with actions of $\Z$) of certain important classes of C*-algebras. Given such an automorphism $\alpha$, we now know from \autoref{thm_main} that $\dr^c(\alpha)$ must be one of $0,1,2$ and $+\infty$. We begin by showing that all these values do indeed occur as the value of $\dr^c(\alpha)$ for some automorphism $\alpha$ of a commutative C*-algebra. \\

We then look at Kirchberg algebras and UHF algebras, where we use the work of Nakamura \cite{nakamura2000} and Kishimoto \cite{kishimoto1995} to show that the value $2$ does not occur for such algebras. Finally, we look at the Jiang--Su algebra, where we explicitly compute both $\dr(\gamma)$ and $\dr^c(\gamma)$ for any automorphism $\gamma \in \Aut(\z)$, using the work of Szab\'{o} \cite{szabo2019a} and Szab\'{o}, Wu and Zacharias \cite{szabo2019}.

\subsection{Commutative C*-Algebras}

For a compact Hausdorff space $X$, every automorphism $\alpha$ of $C(X)$ corresponds to a homeomorphism $T_{\alpha}$ of $X$. Therefore, we may think of the Rokhlin dimension of $\alpha$ as a reflection of the dynamical properties of the pair $(X,T_{\alpha})$. \\

Firstly, if $\dr^c(\alpha) < \infty$, then $T_{\alpha}$ necessarily induces a free action of $\Z$ by \cite[Proposition 4.2]{self_suresh_res_finite}. Therefore, we will now consider those values that could occur if the action is free. \\

Secondly, if $\dr^c(\alpha) = 0$, then $C(X)$ must contain non-trivial projections. So if $X$ is connected, then $\dr^c(\alpha)\geq 1$ must hold. We will now show that for certain kinds of spaces, the value $1$ is also not possible.  \\

For a compact Hausdorff space, let $\pi^1(X) = [X,S^1]$ denote the first cohomotopy group of $X$. It is well known (see, for instance, \cite[Pg. 146]{rordam2000}) that \[
  \pi^1(X) \cong \frac{\U(C(X))}{\U_0(C(X))}
\]
Hence, if $\pi^1(X) = 0$, then every unitary $u \in \U(C(X))$ is of the form $u = e^{2\pi i g}$ for some $g \in C(X,\R)$. This happens when $X$ is a finite, simply connected CW-complex. Indeed, in this case, Bruschlinsky's theorem \cite[Theorem 8.1]{dowker1947} implies that \[\pi^1(X)\cong H^1(X) \cong \Hom(H_1(X), \Z) = 0.\]
The next lemma identifies how this condition affects automorphisms of $C(X)$. 
\begin{lemma}\label{lem_zero_cohomotopy_unitary}
Let $X$ be a compact, connected space with $\pi^1(X) = 0$, and let $\beta \in \Aut(C(X))$ be an automorphism. Then, for any unitary $u \in \U(C(X))$, we have \[
  \|\beta(u) + u\| = 2.
\]
\end{lemma}
\begin{proof}
Let $T : X\to X$ be the homeomorphism corresponding to $\beta$. It is clear that $\|\beta(u) + u\| \leq 2$, so assume $\|\beta(u) + u\| < 2$ so that \[
  \|\beta(u)u^{\ast} + 1\| < 2.
\]
Since $\pi^1(X) = 0$, there is a $g \in C(X,\R)$ such that $u = e^{2\pi i g}$. Then, $\beta(u) = e^{2\pi i (g\circ T)}$ and therefore $h := (g\circ T - g) \in C(X,\R)$ has the property that \[
  \|e^{2\pi i h} + 1\| < 2.
\]
Hence, there is no point $x \in X$ such that $h(x) = 0$. Since $X$ is connected, either $h > 0$ or $h < 0$ must hold. Assume without loss of generality that $h > 0$ holds. Since $X$ is compact, there is a constant $c > 0$ such that \[
  g(T(x)) - g(x) = h(x) \geq c > 0
\]
for all $x \in X$. So for all $n \geq 1$ and all $x \in X$, we have \[
  g(T^n(x)) - g(x) \geq nc.
\]
This contradicts the fact that $g$ must be bounded.
\end{proof}

For our next result, we will be interested in the space $E_{d+1}(\Z/n\Z)$ for some $d\geq 2$. This is a simply connected CW-complex by \cite[Proposition 4.4.3]{matousek2003}, so the previous lemma applies to it.

\begin{theorem}\label{thm_drc_two}
If $\sigma^{(d+1)}$ denotes the natural action of $\Z$ on $C(E_{d+1}(\overline{\Z}))$, then \[
  \dr^c(\sigma^{(d+1)}) = \begin{cases}
    0 &: \text{ if } d = 0, \\
    1 &: \text{ if } d = 1, \text{ and } \\
    2 &: \text{ if } d\geq 2.
  \end{cases}
\]
\end{theorem}
\begin{proof}
Note that $\dr^c(\sigma^{(d+1)}) \leq d$ holds by \autoref{thm_drc_equivalent_definition}. Hence, $\dr^c(\sigma^{(1)}) = 0$. \\

When $d=1$, we appeal to the fact that the topological join construction commutes with taking inverse limits (see \cite[Lemma 1.9]{self_suresh_af}), so that \[
  E_2(\overline{\Z}) \cong \varprojlim E_2(\Z/n!\Z).
\]
Hence, $E_2(\overline{\Z})$ is an inverse limit of connected spaces with surjective connecting maps, and is therefore connected. This implies that $\dr^c(\sigma^{(2)}) \geq 1$, and therefore $\dr^c(\sigma^{(2)}) = 1$. \\

For $d\geq 2$, we know from \autoref{thm_main} that $\dr^c(\sigma^{(d+1)}) \leq 2$, so we need to show that $\dr^c(\sigma^{(d+1)}) \geq 2$. Write $X := E_{d+1}(\overline{\Z})$. Once again, by \cite[Lemma 1.9]{self_suresh_af}, \[
  X \cong \varprojlim E_{d+1}(\Z/n!\Z).
\]
As above, $X$ is connected, so $\dr^c(\sigma^{(d+1)}) \geq 1$. Now suppose that $\dr^c(\sigma^{(d+1)}) \leq 1$, then $\dr^c(\sigma^{(d+1)}, 2\Z) \leq 1$. By \autoref{thm_drc_equivalent_definition}, there would be a $\Z$-equivariant unital $\ast$-homomorphism \[
    \varphi : (C(E_2(\Z/2\Z)),(\sigma^{2\Z})^{(2)})\to (F(C(X)), \sigma^{(d+1)}_{\infty})
  \]
  Now, $E_2(\Z/2\Z) \cong S^1$ and the action of $\Z$ on it corresponds to the antipodal map. If $\zeta \in C(S^1)$ denotes the identity map, then $u := \varphi(\zeta)$ would be a unitary in $F(C(X))$ satisfying \[
    \sigma_{\infty}(u) = -u.
  \]
  Now, $u$ lifts to a sequence $(a_n) \in \ell^{\infty}(\N,C(X))$ with the property that $\lim_{n\to\infty} \|a_n^{\ast}a_n - 1\| = 0$ and $\lim_{n\to\infty} \|a_na_n^{\ast}-1\| = 0$. By \cite[Exercise 2.8]{rordam2000}, there is a sequence of unitaries $(u_n) \in \U(\ell^{\infty}(\N,C(X)))$  such that $\lim_{n\to\infty} \|u_n-a_n\| = 0$. It follows that \[
    \lim_{n\to\infty} \|\sigma^{(d+1)}(u_n) + u_n\| = 0.
  \]
  In particular, there is a unitary $v \in \U(C(X))$ such that $\|\sigma^{(d+1)}(v) + v\| < 2$. Now note that $C(X)$ is an inductive limit of $\Z$-algebras, \[
    (C(X),\sigma^{(d+1)}) \cong \lim_{n\to \infty} (C(E_{d+1}(\Z/n!\Z)), (\sigma^{n!\Z})^{(d+1)}).
  \]
  Once again, by using \cite[Exercise 2.8]{rordam2000}, this implies that there is an $m \in \N$ and a unitary $w \in \U(C(E_{d+1}(\Z/m!\Z)))$ such that \[
    \|(\sigma^{m!\Z})^{(d+1)}(w) + w\| < 2.
  \]
  However, $\pi^1(E_{d+1}(\Z/m!\Z)) = 0$ as mentioned above, so this violates \autoref{lem_zero_cohomotopy_unitary}. We conclude that $\dr^c(\sigma^{(d+1)}) \geq 2$ must hold.
\end{proof}

In particular, each of the values $\{0,1,2,+\infty\}$ is attained for automorphisms on commutative C*-algebras (Note that the identity automorphism of $C(X)$ has infinite Rokhlin dimension).

\subsection{Kirchberg Algebras}

In \cite{nakamura2000}, Nakamura classified aperiodic automorphisms of Kirchberg algebras up to outer conjugacy by their KK-classes. A key step in his proof was to show that such an automorphism necessarily satisfies the \emph{Rokhlin property}. This property had earlier been used extensively by Kishimoto \cite{kishimoto1995}, \cite{kishimoto1996}, Bratteli, Evans and Kishimoto \cite{bratteli1995} and many others. In \cite{hirshberg2015}, Hirshberg, Winter and Zacharias interpreted this property as \emph{Rokhlin dimension zero with two towers}, while also introducing the single-tower definition (and higher dimensional variants). \\

We now give the double-tower definition, so as to be able to use the work of Nakamura to show that for automorphisms of Kirchberg algebras, the only possible values for Rokhlin dimension are $0, 1,$ and $+\infty$. In particular, the value $2$ cannot occur in this context. \\

\begin{defn}
  For any automorphism $\alpha \in \Aut(A)$ of a unital C*-algebra, we say that $\alpha$ has the \emph{Rokhlin property} if, for each finite set $F\subset A, p \in \N$ and $\epsilon > 0$, there are projections $\{e_{0,0}, \ldots, e_{0,p-1}, e_{1,0}, \ldots, e_{1,p}\} \subset A$ such that ~\begin{enumerate}
    \item $\sum_{r=0}^1 \sum_{j=0}^{p+r-1} e_{r,j} = 1$.
    \item $\|\alpha(e_{r,j}) - e_{r,j+1}\| < \epsilon$ for all $0\leq j\leq p+r-2$
    \item $\|[e_{r,j},a]\| < \epsilon$ for all $r,j$ and all $a \in F$.
  \end{enumerate}
\end{defn}

The following result is now contained in \cite[Proposition 2.8]{hirshberg2015}. Indeed, in the proof of that result, the elements $\{g_{j}^{\ell}\}$ constructed in the proof also commute with each other because they belong to the commutative C*-algebra $C^{\ast}(\{e_{r,j}\})$.

\begin{prop}\label{prop_rokhlin_two_towers}\cite[Proposition 2.8]{hirshberg2015} If $\alpha$ has the Rokhlin property, then $\dr^c(\alpha) \leq 1$.
\end{prop}

Let $A$ be a unital C*-algebra, and $\alpha \in \Aut(A)$ be a fixed automorphism. We say that $\alpha$ is \emph{inner} if there is a unitary $u \in \U(A)$ such that $\alpha = \Ad(u)$ (i.e. $\alpha(a) = uau^{\ast}$ for all $a \in A$). We say that $\alpha$ is \emph{outer} if it is not inner, and we say that $\alpha$ is \emph{aperiodic} if $\alpha^k$ is outer for all non-zero $k \in\Z$. \\

It is known that any action with finite Rokhlin dimension is necessarily aperiodic. Nakamura's work not only proves the converse, but also yields the following result.



\begin{prop}\label{cor_kirchberg_automorphism}
Let $A$ be a unital Kirchberg algebra and $\alpha \in \Aut(A)$, then \[
  \dr(\alpha) = \dr^c(\alpha) \in \{0,1,+\infty\}.
\]
\end{prop}
\begin{proof}
Suppose that $\dr(\alpha) < \infty$, then $\alpha^k$ is outer for each $k\neq 0$ by \cite[Proposition 5.2]{self_suresh_res_finite}. In other words, $\alpha$ is aperiodic. By \cite[Theorem 1]{nakamura2000}, $\alpha$ has the Rokhlin property, and hence $\dr^c(\alpha) \leq 1$ by \autoref{prop_rokhlin_two_towers}. \\

Now, $\dr(\alpha) = 0$ if and only if $\dr^c(\alpha) = 0$ (by definition), and therefore the result follows.
\end{proof}

We now determine the exact values of $\dr^c$ for automorphisms of $\O_2$ and $\O_{\infty}$. Recall that two automorphisms $\alpha$ and $\beta$ of a unital C*-algebra $A$ are said to be \emph{outer conjugate} if there is a unitary $u \in \U(A)$ and an automorphism $\delta \in \Aut(A)$ such that $\beta = \Ad(u)\circ \delta\circ \alpha\circ \delta^{-1}$ (see, for instance, \cite{kishimoto1995}). 

\begin{lemma}\label{lem_outer_conjugate}
If $A$ is separable and $\alpha$ and $\beta$ are outer conjugate, then $\dr(\alpha) = \dr(\beta)$ and $\dr^c(\alpha) = \dr^c(\beta)$.
\end{lemma}
\begin{proof}
  We prove the second equality since the proof of the first is identical (using \autoref{defn_dr} for each $H \in \I(\Z)$). Suppose that $d := \dr^c(\alpha) < \infty$, then taking $D = A$ in \autoref{prop_dr_equivalent_definition}, there are $(d+1)$ $\alpha$-equivariant c.c.p. order zero maps \[
  \varphi_0,\varphi_1,\ldots, \varphi_d : (C(\overline{\Z}), \sigma) \to (F(A),\tilde{\alpha}_{\infty})
\]
with pairwise commuting ranges such that $\sum_{\ell=0}^d \varphi_{\ell}(1_{C(\overline{\Z})}) = 1_{F(A)}$. Since $\Ad(u)$ acts trivially on $F(A)$, $\tilde{\beta}_{\infty} = \tilde{\delta}_{\infty}\circ \tilde{\alpha}_{\infty}\circ \tilde{\delta}^{-1}_{\infty}$. It follows that $\dr^c(\beta)\leq d$, witnessed by the maps $\{\psi_0,\psi_1,\ldots, \psi_d\}$, where $\psi_j := \tilde{\delta}_{\infty}\circ \varphi_j$. By symmetry, we conclude that $\dr^c(\alpha) = \dr^c(\beta)$. 
\end{proof}

Before we look at automorphisms on the Cuntz algebras $\mathcal{O}_2$ and $\mathcal{O}_{\infty}$, we need to construct one example of an automorphism on a specific UHF algebra. In what follows, let $\mathcal{Q}$ denote the universal UHF algebra (i.e. $\mathcal{Q}$ is the unique UHF algebra with the property that $K_0(\mathcal{Q})\cong \Q$).

\begin{ex}\label{ex_automorphism_q} We now construct an automorphism $\beta \in \Aut(\mathcal{Q})$ with $\dr(\beta) = 0$ (This will be useful for \autoref{ex_automorphism_cuntz} and for \autoref{thm_uhf_automorphism}). This construction is well-known to the experts, and goes back to the ideas of Kishimoto \cite[Section 5]{kishimoto1995}. \\
  
  For a finite group $\Gamma$, let $\gamma^{\Gamma}$ denote the cyclic shift on $M_{|\Gamma|}$. i.e. $\gamma^{\Gamma} = \Ad(u_{\Gamma})$ where $u_{\Gamma} : \C^{|\Gamma|}\to \C^{|\Gamma|}$ is the map $(x_1,x_2,\ldots, x_{|\Gamma|}) \mapsto (x_{|\Gamma|}, x_1,x_2,\ldots, x_{|\Gamma|-1})$. For each $k \in \N$, let $\Gamma_k := \Z/k!\Z$ and write \[
    \mathcal{Q} \cong \bigotimes_{k=1}^{\infty} M_{|\Gamma_k|} \text{ and } \beta := \bigotimes_{k=1}^{\infty} \gamma^{\Gamma_k}.
  \]
  The eigenvalues of $u_{\Gamma}$ are the $|\Gamma|^{th}$ roots of unity, so it follows from \cite[Section 5]{kishimoto1995} that $\beta$ has the Rokhlin property. In particular, $\dr^c(\beta) \leq 1$. However, it is not hard to see that in fact $\dr^c(\beta) = 0$.\\

   Fix a finite index subgroup $H = n\Z < \Z$. To show that $\dr^c(\beta,H) = 0$, we appeal to the local definition of Rokhlin dimension, as described in \cite[Lemma 2.3]{self_suresh_res_finite} (see also \cite[Proposition 5.5]{szabo2019}). Fix a finite set $F\subset \mathcal{Q}$ and $\epsilon > 0$. For each $k \in \N$, write \[
    A_k := \bigotimes_{i=1}^k M_{|\Gamma_i|}.
  \]
  Then, we may choose $k \in \N$ large enough so that for each $a \in F$, there exists $a' \in A_k$ such that $\|a-a'\| < \epsilon/2$. Choose $K \geq (k+1)$ so that $n \mid K!$. Now define $\varphi : C(\Z/H)\to A_{K}$ by \[
    \varphi(\delta_{\overline{s}}) = \sum_{j\equiv s\mod(n)} e_{j,j}^{(K)},
  \]
  where $e_{j,j}^{(K)} \in M_{|\Gamma_{K}|}$ is the standard $j^{th}$ rank one projection. Note that this is well-defined because $n\mid K!$. It is then easy to see that $\varphi$ is a unital $\ast$-homomorphism and that $\varphi$ is $\beta$-equivariant. Finally, $\varphi$ takes values in $1_{A_{K-1}}\otimes M_{|\Gamma_K|}$, so the image of $\varphi$ commutes with elements of $A_k$, and hence commutes up to an $\epsilon$-error with elements of $F$. It follows from \cite[Lemma 2.3]{self_suresh_res_finite} that $\dr^c(\beta,H) = 0$. This is true for each $H \in \I(\Z)$, so $\dr^c(\beta) = 0$.
\end{ex}

\begin{ex}\label{ex_automorphism_cuntz} For any unital Kirchberg algebra $A$ and any automorphism $\gamma \in \Aut(A)$, we have observed that \[
  \dr(\gamma) < \infty \Leftrightarrow \gamma \text{ is aperiodic} \Leftrightarrow \dr^c(\gamma)\leq 1.
\]
We now identify the exact values of $\dr^c$ for automorphisms of both $\mathcal{O}_2$ and $\mathcal{O}_{\infty}$. 
~\begin{enumerate}
  \item Let $\mathcal{Q}$ denote the universal UHF algebra. In \autoref{ex_automorphism_q}, we constructed an automorphism $\beta \in \Aut(\mathcal{Q})$ with $\dr^c(\beta) = 0$. Then, $\alpha := \id_{\mathcal{O}_2}\otimes \beta$ is an automorphism on $\mathcal{O}_2\otimes \mathcal{Q} \cong \O_2$. By \cite[Theorem 3.2]{self_suresh_res_finite}, $\dr^c(\alpha) \leq \dr^c(\beta)$, so $\dr^c(\alpha) = 0$. Hence, $\O_2$ admits at least one automorphism with Rokhlin dimension zero. \\
  
  However, if $\gamma \in \Aut(O_2)$ is any other aperiodic automorphism, then $\gamma$ is outer conjugate to $\alpha$ by Nakamura's classification theorem (\cite[Theorem 9]{nakamura2000}) because $KK(\O_2,\O_2) = 0$. It follows by \autoref{lem_outer_conjugate} that $\dr^c(\gamma) = 0$. We conclude that for any automorphism $\gamma \in \Aut(\O_2)$, \[
    \dr^c(\gamma) = \begin{cases}
      0 &: \text{ if $\gamma$ is aperiodic} \\
      +\infty &: \text{ otherwise}.
    \end{cases}
  \]
\item Suppose $\alpha\in \Aut(\O_{\infty})$, then it is well known that $\dr^c(\alpha)$ cannot be zero as the following argument shows: We know that $K_0(\mathcal{O}_{\infty}) \cong \Z$ via an isomorphism that maps $[1]_0$ to $1$. Therefore, $\alpha_{\ast} := K_0(\alpha) = \id_{\Z}$. If $\dr^c(\alpha,2\Z) = 0$ were to hold, then there would exist two projections $\{e_0,e_1\} \subset \O_{\infty}$ such that $e_0 + e_1 = 1$ and \[
  \|\alpha(e_i) - e_{i+1}\| < 1, \qquad i \in \{0,1\}.
\]
Therefore, $\alpha(e_i)$ is unitarily equivalent to $e_{i+1}$. If $x := [e_0]_0 \in K_0(\mathcal{O}_{\infty})$, then \begin{equation*}
\begin{split}
  x &= \alpha_{\ast}(x) = [\alpha(e_0)]_0 = [e_1]_0, \text{ and }  \\
  [1]_0 &= [e_0]_0 + [e_1]_0 = 2x.
\end{split}
\end{equation*}
This equation cannot hold in $\Z$, and so $\dr^c(\alpha,2\Z) \geq 1$. Hence, for any automorphism $\gamma \in \Aut(\O_{\infty})$, \[
  \dr^c(\gamma) = \begin{cases}
    1 &: \text{ if $\gamma$ is aperiodic} \\
    +\infty &: \text{ otherwise}.
  \end{cases}
\]
\end{enumerate}
\end{ex}
\subsection{UHF Algebras}

Let $A$ be a unital C*-algebra that admits a unique trace $\tau$. Let $(\mathcal{H},\pi,\xi)$ be the corresponding GNS representation of $\tau$ so that \[
  \tau(a) = \langle \pi(a)\xi,\xi\rangle
\]
for all $a \in A$. Write $\mathcal{M} := \pi(A)''\subset \b(\mathcal{H})$ for the corresponding von Neumann algebra.  \\

For every automorphism $\alpha \in \Aut(A)$, $\tau\circ \alpha = \tau$ (since the trace is unique), so there is a unique $\ast$-automorphism $\overline{\alpha} \in \Aut(\mathcal{M})$ such that $\overline{\alpha}\circ \pi = \pi\circ \alpha$. By the uniqueness of the extension, it follows that $\overline{\alpha^k} = \overline{\alpha}^k$ for each $k \in \Z$. Recall that an automorphism $\alpha \in \Aut(A)$ is said to be \emph{strongly outer} if $\overline{\alpha}^k$ is outer for each $k\neq 0$ (see \cite[Definition 2.7]{matui2012}). \\

We now consider the case where $A$ is a UHF algebra. In this context, Kishimoto classified strongly outer automorphisms in \cite{kishimoto1995} and \cite{kishimoto1996}. In doing so, he showed that any such automorphism must also have the Rokhlin property. Just as in \autoref{cor_kirchberg_automorphism}, this implies the first part of the next theorem. \\

As for the second part, it is well-known that if an algebra admits an automorphism with Rokhlin dimension zero (with a single tower), then that puts certain K-theoretic constraints on the algebra (see, for instance, \cite[Remark 2.2]{hirshberg2015}). A small extension of that idea shows that the only UHF algebra that admits such an automorphism is the universal UHF algebra $\mathcal{Q}$, which we encountered in \autoref{ex_automorphism_q}.

\begin{theorem}\label{thm_uhf_automorphism}
Let $A$ be a UHF algebra and $\alpha \in \Aut(A)$, then
\[
  \dr(\alpha) = \dr^c(\alpha) \in \{0,1,+\infty\}.
\]
Moreover, if this value is finite, then \[
  \dr(\alpha) = \begin{cases}
    0 &: \text{ if } A \cong \mathcal{Q} \\
    1 &: \text{ otherwise}.
  \end{cases}
\]
\end{theorem}
\begin{proof}
If $\dr(\alpha) < \infty$, then by \cite[Proposition 2.3]{liao2016}, $A\rtimes_{\alpha} \Z$ has a unique trace. By \cite[Theorem 1.3]{kishimoto1995}, this implies that $\alpha$ is strongly outer. Conversely, if $\alpha$ is strongly outer, then $\alpha$ has the Rokhlin property by \cite[Theorem 1.3]{kishimoto1995}. By \autoref{prop_rokhlin_two_towers}, we conclude that \[
    \dr(\alpha) < \infty \Leftrightarrow \alpha \text{ is strongly outer } \Leftrightarrow \dr^c(\alpha)\leq 1.
  \] Since $\dr(\alpha) = 0$ holds if and only if $\dr^c(\alpha) = 0$, the first part of the result follows.\\

  Now suppose $\alpha \in \Aut(A)$ is such that $\dr(\alpha) < \infty$, then we need to show that $\dr(\alpha) = 0$ if and only if $A \cong \mathcal{Q}$. ~\begin{enumerate}[label=(\roman*)]
  \item Suppose $\dr(\alpha) = 0$, then we show that $A\cong \mathcal{Q}$. The argument here is similar to that of $\mathcal{O}_{\infty}$ in \autoref{ex_automorphism_cuntz}. Let $\mathfrak{n}$ be the supernatural number associated to $A$ and let $Q(\mathfrak{n})$ be the subgroup of $\Q$ such that \[
  \tau_{\ast} : K_0(A)\to Q(\mathfrak{n})
\]
is an isomorphism. Now, for any $p \in \N$, there are projections $\{e_0, e_1, \ldots, e_{p-1}\} \subset A$ such that $\sum_{j=0}^{p-1} e_j = 1$ and \[
  \|\alpha(e_j) - e_{j+1}\| < 1.
\]
In particular, $\alpha(e_j)$ is unitarily equivalent to $e_{j+1}$, so \[
  \alpha_{\ast}[e_j]_0 = [\alpha(e_j)]_0 = [e_{j+1}]_0
\]
in $K_0(A)$. However, $\tau\circ \alpha = \tau$, so $\tau_{\ast}([e_{j+1}]_0) = \tau_{\ast}([e_j]_0)$. So in $Q(\mathfrak{n})$ we have \[
  1 = \tau_{\ast}([1_A]_0) = \sum_{j=0}^{p-1} \tau_{\ast}([e_j]_0) = p\tau_{\ast}([e_0]_0).
\]
Hence, $p \mid \mathfrak{n}$. This is true for each $p \in \N$, so $Q(\mathfrak{n}) = \Q$, and therefore $A\cong \mathcal{Q}$ must hold.
\item Conversely, suppose $A = \mathcal{Q}$, then we need to show that $\dr(\alpha) = 0$. By \cite[Theorem 1.4]{kishimoto1995} and the first part of the theorem, any two automorphisms of $\mathcal{Q}$ with finite Rokhlin dimension are outer conjugate to each other. In particular, $\alpha$ is outer conjugate to the automorphism $\beta$ constructed in \autoref{ex_automorphism_q}, so $\dr(\alpha) = \dr(\beta) = 0$ by \autoref{lem_outer_conjugate}.
\end{enumerate}
\end{proof}

\subsection{The Jiang--Su Algebra}

We now turn our attention to the Jiang--Su algebra $\z$, which plays a central role in the classification programme. Group actions on this algebra are, as a result, very well-studied. Our goals here are modest; we wish to compute the Rokhlin dimension of automorphisms of this algebra. To that end, we take inspiration from the work of Sato \cite{sato2010} and Matui--Sato \cite{matui2012}, who made use of the de la Harpe--Skandalis determinant from \cite{delaharpe1984}, which we now introduce.

\begin{rem}\label{rem_determinant} Let $A$ be a unital C*-algebra and $\tau : A\to \C$ be a tracial state on $A$. We extend $\tau$ to the (unnormalized) trace $\tau\otimes \mathrm{Tr}$ on $M_n(A)$, which then induces a map \[
  \tau_{\ast} : K_0(A) \to \R.
\]
We define $T := \R/\tau_{\ast}(K_0(A))$ and write $p : \R\to T$ for the quotient map. Write \[
  \U^0_{\infty}(A) := \lim_{n\to \infty} \U^0_n(A)
\]
via the inclusion maps $u \mapsto \begin{pmatrix}
  u & 0 \\
  0 & 1
\end{pmatrix}$. Fix $I := [0,1]\subset \R$ and a piecewise differentiable function $\xi : I \to \U^0_{n}(A)$, and define \[
  \widetilde{\Delta}_{\tau}(\xi) := \frac{1}{2\pi i}\int_{0}^{1} \tau(\xi'(t)\xi(t)^{-1})dt.
\]
Then, $\widetilde{\Delta}_{\tau}$ has the following properties: ~\begin{enumerate}
  \item If $\xi_1, \xi_2 : I \to \U_{\infty}^0(A)$ are two piecewise differentiable functions, then $\widetilde{\Delta}_{\tau}(\xi_1\xi_2) = \widetilde{\Delta}_{\tau}(\xi_1) + \widetilde{\Delta}_{\tau}(\xi_2)$.
  \item If $\|\xi(t) - 1\| < 1$ for all $t \in I$, then \[
    2\pi i \widetilde{\Delta}_{\tau}(\xi) = \tau(\log(\xi(1))) - \tau(\log(\xi(0))).
  \]
  \item $\widetilde{\Delta}_{\tau}(\xi)$ depends only on the homotopy class of $\xi$ (fixing end-points).
\end{enumerate}
Then, the main results of \cite[Section 1]{delaharpe1984} show that $\widetilde{\Delta}_{\tau}$ induces a group homomorphism \[
  \Delta_{\tau} :\U_{\infty}^0(A) \to T
\]
given by $x \mapsto p(\widetilde{\Delta}_{\tau}(\xi))$, where $\xi : [0,1]\to \U_{\infty}^0(A)$ is any path satisfying $\xi(0) = 1$ and $\xi(1) = x$. The function $\Delta_{\tau}$ is called the \emph{de la Harpe--Skandalis determinant}. \\

By \cite[Proposition 2]{delaharpe1984}, this map has the property that for any $y \in M_{\infty}(A)$ with $y^{\ast}=-y$, \[\Delta_{\tau}(\exp(y)) = p\left( \frac{1}{2\pi i} \tau(y)\right).\]
Moreover, it is clear from the definition of $\widetilde{\Delta}_{\tau}$ that if $\alpha \in \Aut(A)$ is an automorphism satisfying $\tau\circ \alpha = \tau$, then \[
  \Delta_{\tau}\circ \alpha = \Delta_{\tau}.
\]
It is these last two properties of $\Delta_{\tau}$ that will be useful to us.
\end{rem}

The next lemma may be thought of as a non-commutative analogue of \autoref{lem_zero_cohomotopy_unitary}. Observe that $\U(\z)$ is path connected because $K_1(\z) = 0$ and it has stable rank one. Moreover, $\z$ carries a unique trace $\tau$ (which must therefore be invariant under any automorphism) and $\tau_{\ast}(K_0(\z)) = \Z$. We refer the reader to \cite{jiang1999} for proofs of these facts.

\begin{lemma}\label{lem_jiang_su_unitary}
For any $\beta \in \Aut(\z)$ and any unitary $u \in \U(\z)$, \[
  \|\beta(u) + u\| = 2.
\]
\end{lemma}
\begin{proof}
Clearly, $\|\beta(u) + u\| \leq 2$, so suppose $\|\beta(u) + u\| < 2$, then for $w := \beta(u)u^{\ast}$, we have \[
  \|w + 1\| < 2.
\]
Therefore, $1 \notin \sigma(w)$, so there is a continuous function $g : \mathbb{T}\setminus\{1\} \to (0,1)$ such that $z = e^{2\pi i g(z)}$ for all $z \in \sigma(w)$. Now, $a := g(w)$ is positive, invertible, $\|a\| < 1$ and \[
  w = e^{2\pi i a}.
\]
By \autoref{rem_determinant}, \[
  \Delta_{\tau}(w) = p\left(\frac{1}{2\pi i}\tau(2\pi i a)\right) = p(\tau(a)).
\]
However, by the properties of the determinant given above, we have \[
  \Delta_{\tau}(w) = \Delta_{\tau}(\beta(u)u^{-1}) = \Delta_{\tau}(\beta(u)) - \Delta_{\tau}(u) = 0.
\]
It follows that $\tau(a) \in \tau_{\ast}(K_0(\z)) = \Z$. Since $\|a\| < 1$, $\tau(a) = 0$ must hold. \\

However, $\tau$ is a faithful trace and $a$ is positive and invertible, so $\tau(a)$ must be strictly positive. This contradiction proves the result.
\end{proof}

We are now in a position to compute the exact values of Rokhlin dimension for an automorphism of the Jiang--Su algebra, both with and without commuting towers. 

\begin{theorem}\label{thm_jiang_su_drc}
For any automorphism $\alpha \in \Aut(\z)$, 
\[
  (\dr(\alpha),\dr^c(\alpha)) = \begin{cases}
    (1,2) &: \text{ if $\alpha$ is strongly outer} \\
    (+\infty, +\infty) &: \text{ otherwise}.
  \end{cases}
\]
\end{theorem}
\begin{proof}
For any automorphism $\alpha \in \Aut(\z)$, we know from \cite[Corollary 1.2]{liao2016} that $\dr(\alpha) < \infty$ if and only if $\alpha$ induces a strongly outer action of $\Z$ on $\z$. \\

In \cite[Remark 11.13]{szabo2019}, the authors construct an automorphism $\gamma \in \Aut(\z)$ with the property that $\dr(\gamma) = 1$. Therefore, $\gamma$ induces a strongly outer action of $\Z$ on $\z$, and by \cite[Theorem 2.14]{szabo2019a}, $\dr^c(\gamma) \leq 2$ must hold. \\

Suppose that $\dr^c(\gamma)\leq 1$, then the proof of \autoref{thm_drc_two} yields a unitary $u \in \U(\z)$ with the property that $\|\gamma(u) + u\| < 2$. This contradicts \autoref{lem_jiang_su_unitary}, so \[
  \dr(\gamma) = 1 \text{ and } \dr^c(\gamma) = 2.
\]
Moreover, by the main results of \cite{szabo2019a}, there is a unique strongly outer action of $\Z$ on $\z$ up to outer conjugacy. Therefore, any such action must have the same Rokhlin dimensions as $\gamma$ by \autoref{lem_outer_conjugate}. This completes the proof.
\end{proof}
In an earlier draft of the paper, we had proved that $\dr^c(\gamma) \in \{2,+\infty\}$ for the example from \cite[Remark 11.13]{szabo2019}. The author is indebted to G\'{a}bor Szab\'{o} for pointing out that the value of $\dr^c(\gamma)$ given above is indeed $2$, thereby upgrading this result from an example to a theorem. We have chosen to refer to the results of \cite{szabo2019a}, but it should be noted that the fact that $\dr^c(\gamma) < \infty$ could have been obtained from the genericity result of \cite[Theorem 3.4]{hirshberg2015} as well.

\subsection*{AI Declaration} The author used an LLM at exactly one point in the preparation of this article: to identify the role of the inflation map (see \autoref{defn_inflation_map}) and prove \autoref{prop_cohomology_isomorphism}. All the other content of this manuscript is solely due to the author (including the mistakes, if any).


\begin{thebibliography}{MV26a}

\bibitem[BEK95]{bratteli1995}
Ola Bratteli, David~E. Evans, and Akitaka Kishimoto.
``The Rohlin Property For Quasi-Free Automorphisms of the Fermion Algebra''.
In: \emph{Proceedings of the London Mathematical Society} s3-71.3 (Nov.\ 1995),
pp.~675--694. \textsc{issn}: 00246115.
\textsc{doi}: \href{https://doi.org/10.1112/plms/s3-71.3.675}{10.1112/plms/s3-71.3.675}.

\bibitem[Bro82]{brown1982}
Kenneth~S. Brown.
\emph{Cohomology of Groups}. Vol.~87. Graduate Texts in Mathematics.
New York, NY: Springer, 1982. \textsc{isbn}: 978-1-4684-9329-0 978-1-4684-9327-6.
\textsc{doi}: \href{https://doi.org/10.1007/978-1-4684-9327-6}{10.1007/978-1-4684-9327-6}.

\bibitem[Con75]{connes1975}
Alain Connes.
``Outer Conjugacy Classes of Automorphisms of Factors''.
In: \emph{Annales scientifiques de l'\'Ecole normale sup\'erieure} 8.3 (1975),
pp.~383--419. \textsc{issn}: 0012-9593, 1873-2151.
\textsc{doi}: \href{https://doi.org/10.24033/asens.1295}{10.24033/asens.1295}.

\bibitem[Die87]{dieck1987}
Tammo tom Dieck.
\emph{Transformation Groups}. De Gruyter Studies in Mathematics v.~8.
Walter De Gruyter, Berlin, 1987. \textsc{isbn}: 978-3-11-009745-0.

\bibitem[Dow47]{dowker1947}
C.~H. Dowker.
``Mapping Theorems for Non-Compact Spaces''.
In: \emph{American Journal of Mathematics} 69.2 (1947), pp.~200--242.
\textsc{issn}: 0002-9327.
\textsc{doi}: \href{https://doi.org/10.2307/2371848}{10.2307/2371848}.
\textsc{jstor}: \href{https://www.jstor.org/stable/2371848}{2371848}.

\bibitem[dS84]{delaharpe1984}
Pierre de~la Harpe and Georges Skandalis.
``D\'eterminant associ\'e \`a une trace sur une alg\`ebre de Banach''.
In: \emph{Annales de l'Institut Fourier} 34.1 (1984), pp.~241--260.
\textsc{issn}: 1777-5310.
\textsc{doi}: \href{https://doi.org/10.5802/aif.958}{10.5802/aif.958}.

\bibitem[Gar17]{gardella2017}
Eusebio Gardella.
``Rokhlin Dimension for Compact Group Actions''.
In: \emph{Indiana University Mathematics Journal} 66.2 (2017), pp.~659--703.
\textsc{issn}: 0022-2518.
\textsc{doi}: \href{https://doi.org/10.1512/iumj.2017.66.5951}{10.1512/iumj.2017.66.5951}.

\bibitem[GHS21]{gardella2021}
Eusebio Gardella, Ilan Hirshberg, and Luis Santiago.
``Rokhlin Dimension: Duality, Tracial Properties, and Crossed Products''.
In: \emph{Ergodic Theory and Dynamical Systems} 41.2 (Feb.\ 2021), pp.~408--460.
\textsc{issn}: 0143-3857, 1469-4417.
\textsc{doi}: \href{https://doi.org/10.1017/etds.2019.68}{10.1017/etds.2019.68}.

\bibitem[GHV22]{gardella2022}
Eusebio Gardella, Ilan Hirshberg, and Andrea Vaccaro.
``Strongly Outer Actions of Amenable Groups on $\mathcal{Z}$-Stable Nuclear
C$^*$-Algebras''.
In: \emph{Journal de Math\'ematiques Pures et Appliqu\'ees} 162 (June 2022),
pp.~76--123. \textsc{issn}: 00217824.
\textsc{doi}: \href{https://doi.org/10.1016/j.matpur.2022.04.003}{10.1016/j.matpur.2022.04.003}.

\bibitem[Hat02]{hatcher2002}
Allen Hatcher.
\emph{Algebraic Topology}. Cambridge; New York: Cambridge University Press, 2002.
\textsc{isbn}: 978-0-521-79160-1 978-0-521-79540-1.

\bibitem[HO84]{herman1984}
Richard~H. Herman and Adrian Ocneanu.
``Stability for Integer Actions on UHF C$^*$-Algebras''.
In: \emph{Journal of Functional Analysis} 59.1 (Oct.\ 1984), pp.~132--144.
\textsc{issn}: 00221236.
\textsc{doi}: \href{https://doi.org/10.1016/0022-1236(84)90056-9}{10.1016/0022-1236(84)90056-9}.

\bibitem[HP15]{hirshberg2015a}
Ilan Hirshberg and N.~Christopher Phillips.
``Rokhlin Dimension: Obstructions and Permanence Properties''.
In: \emph{Documenta Mathematica} 20 (2015), pp.~199--236.
\textsc{issn}: 1431-0635, 1431-0643.
\textsc{doi}: \href{https://doi.org/10.4171/dm/489}{10.4171/dm/489}.

\bibitem[HP26]{hirshberg2026}
Ilan Hirshberg and N.~Christopher Phillips.
``Values of Rokhlin Dimension for Actions of Compact Groups''.
In: \emph{International Mathematics Research Notices} 2026.17 (Sept.\ 2026),
rnag197. \textsc{issn}: 1073-7928.
\textsc{doi}: \href{https://doi.org/10.1093/imrn/rnag197}{10.1093/imrn/rnag197}.

\bibitem[HWZ15]{hirshberg2015}
Ilan Hirshberg, Wilhelm Winter, and Joachim Zacharias.
``Rokhlin Dimension and C$^*$-Dynamics''.
In: \emph{Communications in Mathematical Physics} 335.2 (Apr.\ 2015),
pp.~637--670. \textsc{issn}: 0010-3616, 1432-0916.
\textsc{doi}: \href{https://doi.org/10.1007/s00220-014-2264-x}{10.1007/s00220-014-2264-x}.

\bibitem[Izu04]{izumi2004}
Masaki Izumi.
``Finite Group Actions on C$^*$-Algebras with the Rohlin Property---II''.
In: \emph{Advances in Mathematics} 184.1 (May 2004), pp.~119--160.
\textsc{issn}: 00018708.
\textsc{doi}: \href{https://doi.org/10.1016/S0001-8708(03)00140-3}{10.1016/S0001-8708(03)00140-3}.

\bibitem[JS99]{jiang1999}
Xinhui Jiang and Hongbing Su.
``On a Simple Unital Projectionless C$^*$-Algebra''.
In: \emph{American Journal of Mathematics} 121.2 (1999), pp.~359--413.
\textsc{issn}: 0002-9327.
\textsc{jstor}: \href{https://www.jstor.org/stable/25098671}{25098671}.

\bibitem[Kis95]{kishimoto1995}
Akitaka Kishimoto.
``The Rohlin Property for Automorphisms of UHF Algebras''.
In: \emph{Journal f\"ur die reine und angewandte Mathematik} 465 (1995),
pp.~183--196.

\bibitem[Kis96]{kishimoto1996}
Akitaka Kishimoto.
``The Rohlin Property for Shifts on UHF Algebras and Automorphisms of Cuntz
Algebras''.
In: \emph{Journal of Functional Analysis} 140.1 (Aug.\ 1996), pp.~100--123.
\textsc{issn}: 00221236.
\textsc{doi}: \href{https://doi.org/10.1006/jfan.1996.0100}{10.1006/jfan.1996.0100}.

\bibitem[Lia16]{liao2016}
Hung-Chang Liao.
``A Rokhlin Type Theorem for Simple C$^*$-Algebras of Finite Nuclear
Dimension''.
In: \emph{Journal of Functional Analysis} 270.10 (May 2016), pp.~3675--3708.
\textsc{issn}: 00221236.
\textsc{doi}: \href{https://doi.org/10.1016/j.jfa.2016.01.014}{10.1016/j.jfa.2016.01.014}.

\bibitem[Mat03]{matousek2003}
Ji\v{r}\'{\i} Matou\v{s}ek.
\emph{Using the Borsuk--Ulam Theorem: Lectures on Topological Methods in
Combinatorics and Geometry}. Universitext. Berlin; New York: Springer, 2003.
\textsc{isbn}: 978-3-540-00362-5.

\bibitem[Mil56]{milnor1956}
John Milnor.
``Construction of Universal Bundles, II''.
In: \emph{The Annals of Mathematics} 63.3 (May 1956), pp.~430--436.
\textsc{issn}: 0003486X.
\textsc{doi}: \href{https://doi.org/10.2307/1970012}{10.2307/1970012}.
\textsc{jstor}: \href{https://www.jstor.org/stable/1970012}{1970012}.

\bibitem[MS12]{matui2012}
Hiroki Matui and Yasuhiko Sato.
``$\mathcal{Z}$-Stability of Crossed Products by Strongly Outer Actions''.
In: \emph{Communications in Mathematical Physics} 314.1 (Aug.\ 2012),
pp.~193--228. \textsc{issn}: 1432-0916.
\textsc{doi}: \href{https://doi.org/10.1007/s00220-011-1392-9}{10.1007/s00220-011-1392-9}.

\bibitem[MV26a]{self_suresh_af}
Sureshkumar M and Prahlad Vaidyanathan.
``Rokhlin Dimension and Inductive Limit Actions on AF-Algebras''.
In: \emph{Annals of Functional Analysis} 17.2 (Mar.\ 2026), p.~35.
\textsc{issn}: 2008-8752.
\textsc{doi}: \href{https://doi.org/10.1007/s43034-026-00502-3}{10.1007/s43034-026-00502-3}.

\bibitem[MV26b]{self_suresh_res_finite}
Sureshkumar M and Prahlad Vaidyanathan.
``Rokhlin Dimension: Permanence Properties and Ideal Separation''.
In: \emph{Groups, Geometry, and Dynamics} 20.2 (2026), pp.~555--588.
\textsc{issn}: 1661-7207.
\textsc{doi}: \href{https://doi.org/10.4171/ggd/825}{10.4171/ggd/825}.

\bibitem[Nak00]{nakamura2000}
Hideki Nakamura.
``Aperiodic Automorphisms of Nuclear Purely Infinite Simple C$^*$-Algebras''.
In: \emph{Ergodic Theory and Dynamical Systems} 20.6 (Dec.\ 2000),
pp.~1749--1765. \textsc{issn}: 0143-3857, 1469-4417.
\textsc{doi}: \href{https://doi.org/10.1017/S0143385700000973}{10.1017/S0143385700000973}.

\bibitem[RLL00]{rordam2000}
M.~R{\o}rdam, F.~Larsen, and N.~Laustsen.
\emph{An Introduction to K-Theory for C$^*$-Algebras}. Vol.~49. London
Mathematical Society Student Texts. Cambridge University Press, 2000.
\textsc{doi}: \href{https://doi.org/10.1017/CBO9780511623806}{10.1017/CBO9780511623806}.

\bibitem[RZ10]{ribes2010}
Luis Ribes and Pavel Zalesskii.
\emph{Profinite Groups}. Berlin, Heidelberg: Springer Berlin Heidelberg, 2010.
\textsc{isbn}: 978-3-642-01641-7 978-3-642-01642-4.
\textsc{doi}: \href{https://doi.org/10.1007/978-3-642-01642-4}{10.1007/978-3-642-01642-4}.

\bibitem[Sat10]{sato2010}
Yasuhiko Sato.
``The Rohlin Property for Automorphisms of the Jiang--Su Algebra''.
In: \emph{Journal of Functional Analysis} 259.2 (July 2010), pp.~453--476.
\textsc{issn}: 00221236.
\textsc{doi}: \href{https://doi.org/10.1016/j.jfa.2010.04.006}{10.1016/j.jfa.2010.04.006}.

\bibitem[SWZ19]{szabo2019}
G\'abor Szab\'o, Jianchao Wu, and Joachim Zacharias.
``Rokhlin Dimension for Actions of Residually Finite Groups''.
In: \emph{Ergodic Theory and Dynamical Systems} 39.8 (Aug.\ 2019),
pp.~2248--2304. \textsc{issn}: 0143-3857, 1469-4417.
\textsc{doi}: \href{https://doi.org/10.1017/etds.2017.113}{10.1017/etds.2017.113}.

\bibitem[Sza19]{szabo2019a}
G\'abor Szab\'o.
``Actions of Certain Torsion-Free Elementary Amenable Groups on Strongly
Self-Absorbing C$^*$-Algebras''.
In: \emph{Communications in Mathematical Physics} 371.1 (Oct.\ 2019),
pp.~267--284. \textsc{issn}: 0010-3616, 1432-0916.
\textsc{doi}: \href{https://doi.org/10.1007/s00220-019-03435-2}{10.1007/s00220-019-03435-2}.

\end{thebibliography}
\end{document}